\documentclass[12pt]{article}

\usepackage{amsfonts}
\usepackage{amsmath}
\usepackage{amssymb, latexsym}
\usepackage{graphicx}
\usepackage[usenames]{color}

\usepackage{mathrsfs}

\newtheorem{theorem}{Theorem}[section]
\newtheorem{proposition}[theorem]{Proposition}

\newtheorem{lemma}[theorem]{Lemma}
\newtheorem{remark}[theorem]{Remark}
\newtheorem{definition}[theorem]{Definition}

\numberwithin{equation}{section} 

\def\Sc{{Schr\"odinger} }
\def\cala{{\mathcal{A}}}
\def\calf{{\mathcal{F}}}

\def\calh{\mathcal{H}}
\def\calx{\mathcal{X}}
\def\calb{{\mathcal{B}}}
\def\calg{{\mathcal{G}}}
\def\calu{{\mathcal{U}}}
\def\bbc{{\mathbb{C}}}
\def\bbe{{\mathbb{E}}}
\def\bbr{{\mathbb{R}}}
\def\bbp{{\mathbb{P}}}

\def\en t{{{\rm Z}\mkern-5.5mu{\rm Z}}}

\def\oo{{\omega}}

\def\<{\left<}
\def\>{\right>}
\def\({\left(}
\def\){\right)}

\def\D{\Delta }
\def\9{{\infty}}
\def\barr{\begin{array}}
\def\earr{\end{array}}
\def\ov{\overline}
\def\wt{\widetilde}
\def\wh{\widehat}
\def\ol{\overline}

\def\vf{{\varphi}}
\def\lbb{{\lambda}}

\def\a{{\alpha}}

\def\n{\noindent }

\def\D{{\Delta}}

\def\3{\subset }
\def\na{{\nabla}}
\def\sk{\smallskip }

\def\bk{\bigskip }

\def\ve{{\varepsilon}}

\begin{document}

\begin{center}
{\Large{\bf The stochastic logarithmic \Sc equation}}
\bigskip\bk

{\large{\bf Viorel Barbu}}\footnote{Octav Mayer Institute of
Mathematics (Romanian Academy)   and Al.I. Cuza University and,
700506, Ia\c si, Romania. This work was supported by the DFG through
CRC 701 and by CNCS-VEFISCDI (Romania) project PN-II-2012-4-0456.},
{\large{\bf Michael R\"ockner}}\footnote{Fakult\"at f\"ur
Mathematik, Universit\"at Bielefeld,  D-33501 Bielefeld, Germany.
This research was supported by the DFG through CRC 701.},
{\large{\bf Deng Zhang}}\footnote{Department of Mathematics,
Shanghai Jiao Tong University, 200240 Shanghai, China. }
\end{center}

\bk\bk\bk

\begin{quote}
\n{\small{\bf Abstract.} In this paper we prove global existence and
uniqueness of solutions to the stochastic logarithmic Schr\"odinger
equation with linear multiplicative noise. Our approach is mainly
based on the rescaling approach and the method of maximal monotone
operators. In addition, uniform estimates of solutions in
the energy space $H^1(\bbr^d)$ and in an appropriate Orlicz space are also obtained here.} \\


\end{quote}

\vfill

\section{Introduction and main result.}

The logarithmic Schr\"{o}dinger equation
\begin{align} \label{deter-equa-x}
    i\frac{d u}{dt} + \D u +  u \log |u|^2 = 0, \ \ in\ \bbr^+ \times
    \bbr^d,
\end{align}
has wide applications in quantum mechanics, quantum optics, nuclear
physics, open quantum systems, Bose-Einstein condensation and so on.
It was first proposed in \cite{BM76} as a model of nonlinear wave
mechanics.  As a matter of fact, as shown in \cite{BM76}, the
logarithmic nonlinearity arising in \eqref{deter-equa-x} is the
unique nonlinearity for which the separability hypothesis of
noninteracting subsystems of the Schr\"odinger theory holds. It also
possesses many other attractive features, including the additivity
of the energy for noninteracting subsystems, the validity of the
lower energy bound and Planck's relation for all stationary states.
All these make this equation unique among nonlinear wave equations.
See e.g. \cite{BM76, BM79, GLN10, Z10}.  We also refer to \cite{L80,
N85} for the derivation of this equation from Nelson's stochastic
quantum mechanics \cite{N66}.

Motivated by the physical significance above, we are here mainly
concerned with well-posedness of the logarithmic Schr\"odinger
equation in the stochastic case, that is,
\begin{align} \label{equa-x}
   idX&=\Delta Xdt   + \lbb X\log|X|^2dt  -i\mu Xdt + iXdW,\ t\in(0,T), \nonumber \\
   X(0)&=x \in L^2,
\end{align}
Here,  $\lbb\in \mathbb{R}$,  $W$ is the Wiener process
\begin{align} \label{W}
   W(t,\xi) = \sum\limits_{j=1}^n \mu_j e_j(\xi) \beta_j(t),\ t\geq
   0,\ \xi\in\mathbb{R}^d,
\end{align}
where $d\geq 1$, $\{\mu_j\}_{j=1}^n$ are complex numbers,
$\{e_j\}_{j=1}^n$ are real-valued functions, and
$\{\beta_j\}_{j=1}^n$ is a family of independent real valued
Brownian motions on a probability space $(\Omega, \mathscr{F},
\mathbb{P})$ with normal (in particular right-continuous) filtration
$(\mathscr{F}_t)_{t\geq 0}$. For simplicity, we assume that $n<\9$.

Moreover,
\begin{align} \label{mu}
   \mu(\xi) = \frac12 \sum\limits_{j=1}^n |\mu_j|^2 e_j^2(\xi),\ \
   \xi\in\mathbb{R}^d.
\end{align}

The stochastic equation \eqref{equa-x} can be derived from
\eqref{deter-equa-x} with an additional potential $VX$, where the
random potential $V$ fluctuates rapidly and so can be approximated
by the Gaussian noise $\dot{W}$. Moreover, the linear multiplicative
noise $iXdW$ together with the term $-i\mu X dt$ also plays an
important role in the theory of measurements continuous in time in
open quantum systems. In this case, one main feature is that
$|X(t)|_2^2$ is a continuous martingale. This fact implies the mean
norm square conservation of $X(t)$ and allows to define a new
probability law, the ``physical'' probability law, which has
important applications to open quantum systems. For more physical
interpretations, we refer to \cite{BG09}, \cite{BRZ14, BRZ14.2} and
the references therein.

The stochastic nonlinear Schr\"{o}dinger equation with the
polynomial nonlinearity $\lbb|X|^{\a-1}X$ was first studied in
\cite{BD99, BD03}, based on the mild formulation of the stochastic
equation. The optimal exponents of the nonlinearity for the global
well-posedness were recently achieved in \cite{BRZ14, BRZ14.2},
based on the rescaling transformation (see \eqref{rescal} below) and
the Strichartz estimates established in \cite{MMT08} for lower order
perturbations of the Laplacian. However, the contraction mapping
arguments used in the mentioned works are not applicable here, due
to the fact that the function $y\rightarrow y\log|y|^2$ is not
locally Lipschitz.

One of the main features of the logarithmic nonlinearity is the
quasi-monotonicity. Based on this, the global well-posedness of the
deterministic equation \eqref{deter-equa-x} was first studied in
\cite{CH80} in the distribution sense for initial data in $L^2$ or
$H^1$. Later, the global well-posedness was also proved in
\cite{C83} for initial data in $H^1$ and in some convenient Orlicz
space, which is closely related to the logarithmic nonlinearity. We
also refer to \cite{GLN10} for the global well-posedness for initial
data in $H^1$ with finite momentum.

Furthermore, stochastic partial differential equations with monotone
coefficients are also extensively studied in the literature. We
refer to \cite{KR79}, \cite{PR07}, \cite{RRW07} and the references
therein. Recently, based on the rescaling approach and operatorial
reformulation, the approach of maximal monotone operators was
developed in \cite{BR14} in a general infinite dimensional setting,
which has applications to new existence and uniqueness results of
various stochastic models with linear multiplicative noise.

Inspired by the quasi-monotone feature of the logarithmic
nonlinearity and the works mentioned above, we shall employ the
rescaling transformation and the method of maximal monotone
operators to study the global well-posedness of \eqref{equa-x}.

However, it should be mentioned that, the results in \cite{BR14} are
not applicable here, since the operator $i\Delta$ in \eqref{equa-x}
is not coercive (see \cite[(2.3)]{BR14}).

Moreover, another difficulty arises from the passage to the limit in
the approximating equation (see \eqref{app-equa-y} below). Because
even if a space $\mathcal{X}$ is compactly imbedded into another one
$\mathcal{Y}$, we generally do no have the compact imbedding from
$L^p(\Omega; \mathcal{X})$ to $L^p(\Omega; \mathcal{Y})$, $1\leq p
\leq \9$, the classical deterministic method as in
\cite{C83,C03,CH80,GLN10} to pass to the limit in the nonlinear term
can not directly be applied here.

In order to overcome these difficulties, inspired by \cite{C83,
C03}, we will consider the initial data in the energy space
$H^1(\bbr^d)$ and an appropriate Orlicz space $V$ (see \eqref{def-V}
below). These spaces allow to control the singularity of the
logarithmic nonlinearity at infinity and at the origin respectively.
More importantly, they are also suitable spaces for the maximal
monotonicity of the logarithmic nonlinearity, which makes the
passage to the limit in the approximating equation possible, thereby
yielding the global well-posedness.

To state our results precisely, let us first introduce some
necessary notations. Take $H=L^2(\mathbb{R}^d; \mathbb{C})=:L^2$
with the scalar product defined by $   \<u,v\> = \int_{\mathbb{R}^d}
u\overline{v} d\xi, \ u,v\in H$, and the norm $|u|_2=\<u,u\>^{\frac
12}$. Let $H^1$ denote the classical Sobolev space, i.e. $H^1=\{u\in
L^2: \na u\in L^2\}$ with norm $|u|^2_{H^1} = |u|_2^2 + |\na
u|^2_2$. We also use the standard notation $L^p= L^p(\bbr^d)$,
$1\leq p \leq \9$, for the space of all $p$-integrable complex
functions with the norm $|\cdot|_{L^p}$.

Moreover, as in \cite{C83}, define the function
\begin{align} \label{def-N}
   N(x) = \left\{
            \begin{array}{ll}
              -x^2 \log x^2, & \hbox{if $0\leq x\leq e^{-3}$;} \\
              3x^2 + 4e^{-3}x - e^{-6}, & \hbox{if $e^{-3}\leq x$.}
            \end{array}
          \right.
\end{align}
$N$ is a positive convex and increasing function, and $N\in
C^1([0,\9)) \cap C^2 ((0,\9))$. The Orlicz space $V$ corresponding
to $N$ is defined by
\begin{align} \label{def-V}
    V=\{ u\in L^1_{loc}: N(|u|)\in L^1 \},
\end{align}
equipped with the Luxembourg norm
\begin{align}
    \|u\|_V = \inf \{ k>0: \int N(k^{-1}|u(\xi)|)d\xi \leq 1  \}.
\end{align}
Here as usual $L^1_{loc}$ is the space of all locally Lebesgue
integrable functions. It is proved in \cite[Lemma 2.1]{C83} that $N$
is a Young-function which is $\bigtriangleup_2$-regular and  $(V,
\|\cdot\|_{V})$ is a separable reflexive Banach space (see also
\cite{C03} and \cite{A75}). We also have that (see
\cite[$(2.2)$]{C83}) for any $u\in V$,
\begin{align} \label{V-N}
   \min\{ \|u\|_V,  \|u\|^2_V \} \leq \int N(|u(\xi)|) d\xi
   \leq \max \{ \|u\|_V,  \|u\|^2_V \}
\end{align}

Now, set $U:= H^1 \cap V$. $U$ is a reflexive Banach space equipped
with the norm $\|u\|_{U}= |u|_{H^1} + \|u\|_V$, for any $u \in U$,
and its dual space is $U'=H^{-1}+ V'$ with the norm $\|u\|_{U'}=\inf
\{ |u_1|_{H^{-1}} + \|u_2\|_{V'}: u= u_1+u_2, u_1\in H^{-1}, u_2\in
V' \}$. One advantage for introducing the space $U$ is that the
nonlinear operator $u \mapsto u \log |u|^2$ is continuous from $U$
to $U'$ (see \cite[Lemma 2.6]{C83}).\\

The precise definition of solutions to \eqref{equa-x} is given
below.

\begin{definition} \label{Def-x}
A continuous $H$-valued $(\calf_t)$-adapted process $X$ is said to
be a solution to \eqref{equa-x} if for any $p\geq 3$, $X\in
L^p(\Omega \times(0,T); U)$, $X\log |X|^2 \in L^{p'}(\Omega \times
(0,T); U')$, and it satisfies $\bbp$-a.s. for all $t\in[0,T]$
\begin{align} \label{equa-x*}
    X(t)=& x - \int_0^t \(i\D X(s) ds + \mu X(s)
           + \lbb X(s)\log| X(s)|^2 \) ds
          +  \int_0^t X(s)dW(s),
\end{align}
where the stochastic term is taken in It\^o's sense.
\end{definition}

We also assume that the spatial functions $\{e_j\}_{j=1}^n$ in the
noise $W$ satisfy the hypothesis: {\it
\begin{itemize}
  \item[{\rm(H)}] $e_j\in \calb^\9(\bbr^d)$ such that for each $1\leq k \leq
  d$, $1\leq m \leq n$,
   $$ |\partial_k e_m(\xi)| \leq \lbb(|\xi|),\ \ \xi\in \bbr^d ,$$
  where $\calb^\9 = \{ f\in C^\9(\bbr^d),
\partial_\alpha f\in L^\9,\ for\ all\ \alpha\}$, and $\lbb(\cdot) $ is a positive non-increasing function in $C([0,\9]) \cap L^1
  ([0,\9))$.
\end{itemize}}

The main result of this article is formulated as follows.
\begin{theorem} \label{thm-x}
Under Hypothesis $(H)$, for any initial datum $x\in U$ and $0<T<\9$,
there exists a unique solution $X$ to \eqref{equa-x} in the sense of
Definition \ref{Def-x}.

Moreover, for any $p\geq 2$,
\begin{align}
\label{integ-l2-x}
    &\bbe \|X(t)\|^{p}_{L^\9(0,T;U)} <\9,
\end{align}
\begin{align}
    & \bbe \| X(t) \log |X(t)|^2
    \|^p_{L^\9(0,T;U')}<\9,
\end{align}
and
\begin{align}
    \bbe  \| e^{W(t)} \frac{d}{dt} (e^{-W(t)}X(t))
    \|^p_{L^\9(0,T;U')} < \9.
\end{align}
\end{theorem}

The remainder of this paper is organized as follows. In Section
\ref{Oper-Equa}, we first apply the rescaling transformation to
reduce the original stochastic equation \eqref{equa-x} to a random
equation (see \eqref{equa-y}), and then we introduce some
appropriate spaces and prove the maximal monotonicity of the
logarithmic nonlinearity. Section \ref{App-Equa} is mainly concerned
with the approximating equation. We first obtain the $H^1$-global
well-posedness and derive the uniform estimate in the energy space
in Subsection \ref{H1-GWP}. Then in Subsection \ref{Uni-Log}, in
order to control the singularity of the logarithmic nonlinearity at
the origin, we start with the analysis of the entropy function and
then prove the uniform estimates in the Orlicz space. Section
\ref{Proof-Main} is mainly devoted to the proof of the main result.
As mentioned above, the maximal monotonicity will play an important
role in the passage to the limit in the approximating equation. Some
technical details are postponed to the Appendix.

Throughout this paper, $C$ denotes various constants which may
change from line to line.

\section{Random equation} \label{Oper-Equa}

Taking into account the quasi-monotone feature of the logarithmic
nonlinearity, we first use the change of variable $X\to e^{-2|\lbb|
t}X$ to reformulate the original equation \eqref{equa-x} as
\begin{align} \label{equa-x'}
    idX&=\Delta Xdt   + \lbb X\log|X|^2dt  +(4\lbb|\lbb|t -2i|\lbb| -i\mu) Xdt + iXdW, \nonumber \\
   X(0)&=x \in L^2.
\end{align}
Then, applying the rescaling transformation
\begin{align} \label{rescal}
   X=e^W y,
\end{align}
which can be seen as a Doss-Sussman transformation generalized to
infinite dimensions, we can reduce the stochastic equation
\eqref{equa-x'} to a random Schr\"odinger equation
\begin{align} \label{equa-y}
   &\frac{ d y}{d t}(t) = - i e^{-W(t)}\D (e^{W(t)}y(t))
    - ( 2|\lbb| + 4i\lbb|\lbb|t + \wh{\mu} ) y(t) \nonumber \\
    &\qquad \qquad - \lbb i y(t) \log|e^{W(t)} y(t)|^2,\ \ a.e.\ t\in(0,T),  \\
   &y(0)=x, \nonumber
\end{align}
where $\wh \mu = \frac 12 \sum\limits_{j=1}^N (\mu_j^2 + |\mu_j|^2)
e_j^2.$

In order to formulate the definition of solutions to \eqref{equa-y},
proceeding as in \cite{BR14}, we consider the Hilbert space
$\mathcal{H}$ of all $H$-valued $(\mathscr{F}_t)_{t\geq 0}$-adapted
processes $y:[0,T] \rightarrow H$ with the scalar product
\begin{align} \label{def-h}
   \<y,z\>_{\mathcal{H}} = \mathbb{E} \int_0^T
   \<e^{W(t)}y(t),e^{W(t)}z(t)\>dt,
\end{align}
and the norm
\begin{align*}
  |y|_{\calh} = \(\mathbb{E}\int_0^T|e^{W(t)}y(t)|_H^2dt\)^{\frac 12}.
\end{align*}

For any $p\geq 3$, consider the space $\calu$ of all
$(\mathscr{F}_t)_{t\geq 0}$-adapted processes $y:= [0,T] \to U$ such
that
\begin{align} \label{def-calU}
    \| y \|_{\calu}^p = \bbe \int_0^T \| e^{W(t)} y(t) \|^{p}_U dt
    <\9,
\end{align}
Let $\calu'$ denote the dual space of $\calu$. In fact, $\calu'$ is
the space of all $(\mathscr{F}_t)_{t\geq 0}$-adapted processes
$y:[0,T] \to \calu'$ such that
\begin{align}  \label{def-calU'}
    \| y \|_{\calu'}^{p'} = \bbe \int_0^T \| e^{W(t)} y(t) \|^{p'}_{U'}
    dt
    <\9.
\end{align}
We have $\calu \subset \calh \subset \calu'$, algebraically and
topologically.

Set
\begin{align} \label{Def-G}
 (\calg y)(t):=& \lbb i y(t) \log|e^{W(t)}y(t)|^2 + 2|\lbb|y(t),\ \  y\in
 D(\calg)=\mathcal{U},
\end{align}

Analogously to Definition \ref{Def-x}, the solutions to
\eqref{equa-y} is now defined below.

\begin{definition} \label{def-y}
A solution to \eqref{equa-y} is a continuous $H$-valued
$(\mathscr{F}_t)$-adapted progress $y$, such that $y\in \calu$,
$y\log|e^Wy|^2 \in \calu'$, and it satisfies $\bbp$-a.s. for all
$t\in[0,T]$
\begin{align} \label{equa-y*}
    y(t) = x - \int_0^t \( i e^{-W(s)}\D (e^{W(s)} y(s))
           + (4i\lbb|\lbb|t + \wh{\mu}) y(s) + \calg(y(s))\) ds.
\end{align}
\end{definition}

We refer to \cite[Lemma 8.1]{BR14} for a rigorous proof of the
equivalence of solutions to \eqref{equa-x*} and \eqref{equa-y*}.
Therefore, the proof of Theorem \ref{thm-x} is now reduce to the
theorem as follows.

\begin{theorem}\label{thm-y}
Under Hypothesis $(H)$, for any initial datum $x\in U$ and $0<T<\9$,
there exists a unique solution $y$ to \eqref{equa-y*} in the sense
of Definition \ref{def-y}.

Moreover, for all $p\geq 2$,
\begin{align} \label{Bdd-U-y}
    &\bbe \|e^{W(t)}y(t)\|^{p}_{L^\9(0,T;U)} <\9,
\end{align}
\begin{align} \label{Bdd-U'-log}
    &\bbe  \| e^{W(t)}y(t) \log |e^{W(t)}y(t)|^2
    \|^{p}_{L^\9(0,T;U')} <\9,
\end{align}
and
\begin{align} \label{Bdd-U'-dy}
     \bbe
     \|e^{W(t)}\frac{d}{dt}y(t)\|^p_{L^\9(0,T;U')} <\9.
\end{align}
\end{theorem}

The remainder of this paper is devoted to the proof of Theorem
\ref{thm-y}. We will mainly consider the case
$d\geq 3$. The simpler cases $d=1,2$ can be proved similarly.\\

In the end of this section, let us show the maximal monotonicity of
the operator $ \calg$.  Recall that an operator $A:\calx\to\calx'$
(possibly nonlinear) from a Banach space $\calx$ to its dual
$\calx'$ is said to be  monotone if
\begin{align*}
   Re\ {}_{\calx'}\<Ay_1-Ay_2,y_1-y_2\>_{\calx} \geq 0,\ \forall y_1,y_2\in D(A),
\end{align*}
and maximal monotone if it has no nontrivial monotone extensions in
$\calx \times \calx'$.

\begin{proposition} \label{Max-ABG}
For any $p\geq 3$, the operator $\calg$ is maximal monotone from
$\calu$ to $\calu'$.
\end{proposition}

{\it \bf Proof.} In view of \cite[Theorem $2.4$]{B10} the
maximality, since the demicontinuity implies the hemicontinuity, it
suffices to prove that $\calg$ is monotone and demicontinous from
$\calu$ to $\calu'$, i.e., if $y_n, y\in \calu$ such that $y_n\to y$
in $\calu$, then
\begin{align} \label{hemi-g}
    {}_{\calu'} \<\calg(y_n), z \>_{\calu} \to  {}_{\calu'} \<\calg(y),
    z\>_{\calu},\ \ z\in \calu.
\end{align}

For this purpose, we first note that by the definition of $\calg$ in
\eqref{Def-G},
\begin{align*}
   &Re\ {}_{\calu'} \< \calg(y_1)-\calg(y_2), y_1-y_2\>_{\calu} \\
   =&2|\lbb||y_1-y_2|_{\calh}^2 -2 \lbb Im\ {}_{\calu'} \< y_1 \log|e^{W}y_1| - y_2 \log|e^{W}y_2|, y_1-y_2
   \>_{\calu} \geq 0.
\end{align*}
where in the last step we used \eqref{mono-lve-0} below with $\ve
=0$, and so the monotonicity of $\calg$ follows.

In order to prove the demicontinuity  \eqref{hemi-g}, we will show
that
\begin{align} \label{bdd-lp'-u'}
   \| e^W \calg(y_n)\|_{L^{p'}(\Omega \times (0,T);U')} \leq C < \9,
\end{align}
where $C$ is independent of $n$. Then, for any subsequence of $\{n\}
\to \9$, there exists a further subsequence (still denoted by
$\{n\}$) such that $e^W \calg(y_n) \overset{\omega}{\rightharpoonup}
\eta$, in $L^{p'}(\Omega \times (0,T); U')$, where $
\overset{\omega}{\rightharpoonup}$ stands for weak convergence. But,
since $y_n \to y$ in $\calu$, we have $e^W\calg(y_n) \to e^W
\calg(y)$ in measure $\bbp \otimes dt \otimes d\xi$. Hence, we
conclude that $\eta=e^W \calg(y)$, which implies \eqref{hemi-g},
since the subsequence was arbitrary.

It remains to prove \eqref{bdd-lp'-u'}. Set $X_n:= e^W y_n$ and
$L(|X_n|^2):= \log|X_n|^2$.  By the definition of $U'$ and $\calg$
we have
\begin{align} \label{split-u'}
     &\| e^W \calg(y_n)\|_{L^{p'}(\Omega \times (0,T);U')} \nonumber  \\
     \leq& 2|\lbb|\ \|  X_n \|_{L^{p'}(\Omega \times (0,T);H^{-1})}
          +|\lbb| \| I_{\{|X_n|  > e^{-3}\}} X_n L(|X_n|^2) \|_{L^{p'}(\Omega \times (0,T);H^{-1})}  \nonumber \\
     &     + |\lbb| \|I_{\{|X_n| \leq e^{-3}\}}  X_n L(|X_n|^2)   \|_{L^{p'}(\Omega \times (0,T);V')}.
\end{align}

Since for each  $\xi \in \{|X_n| > e^{-3} \}$, $
     |X_n (\xi) L(|X_n(\xi)|^2)| \leq C_\delta(|X_n(\xi)| + |X_n(\xi)|^{1+\delta})$ with $C_\delta$
independent of $n$. By Sobolev's imbedding theorem with $\delta \leq
\frac{2}{d-2}$,
\begin{align*}
    & \|  X_n \|_{L^{p'}(\Omega \times (0,T);H^{-1})}
       + \| I_{\{|X_n|  > e^{-3}\}} X_n L(|X_n|^2) \|_{L^{p'}(\Omega \times (0,T);H^{-1})} \nonumber \\
    \leq& C (\|X_n\|_{L^{p'}(\Omega\times(0,T); L^2)}
          + \|X_n\|^{1+\delta}_{L^{(1+\delta)p'}(\Omega\times(0,T);
          L^{2(1+\delta)})}) \\
    \leq& C (\|X_n\|_{L^{p'}(\Omega\times(0,T); L^2)}
          +   \|X_n\|^{1+\delta}_{L^{(1+\delta)p'}(\Omega\times(0,T);
          H^1)}).
\end{align*}
Then, taking $\delta$ such that $0<\delta<p-2$, we have
$(1+\delta)p' < p$ and, via the H\"{o}lder inequality,
\begin{align} \label{bdd-lp'-h-1}
    &    \|  X_n \|_{L^{p'}(\Omega \times (0,T);H^{-1})}
        +\| I_{\{|X_n|  > e^{-3}\}} X_n L(|X_n|^2) \|_{L^{p'}(\Omega \times (0,T);H^{-1})}
         \nonumber \\
    \leq&  C_T ( \| X_n \|_{L^{p}(\Omega \times
          (0,T);L^2)}
          +  \| X_n \|^{1+\delta}_{L^{p}(\Omega \times
          (0,T);H^1)} )
    \leq C_T<\9,
\end{align}
where $C_T$ is independent of $n$.

On the other hand, for each $\xi \in \{|X_n| \leq e^{-3} \}$, as in
the proof of \cite[Lemma 2.5]{C83} we have
\begin{align} \label{N'-N}
\wt{N}(|X_n(\xi) L(|X_n(\xi)|^2 )|) \leq 2 N(|X_n(\xi)|),
\end{align} where $\wt{N}$ is the convex conjugate of $N$.
Then, since $\wt{N}(0)=0$, by \eqref{V-N},
\begin{align} \label{wtN-V}
    \int \wt{N}(I_{\{|X_n|\leq e^{-3}\}} |X_n L(|X_n|^2)|) d\xi
   =&\int  I_{\{|X_n|\leq e^{-3}\}} \wt{N}( |X_n L(|X_n|^2)|) d\xi \nonumber \\
   \leq& 2 \int  I_{\{|X_n|\leq e^{-3}\}} N(|X_n|) d\xi \nonumber \\
   \leq& 2 \max\{\|X_n\|_V, \|X_n\|^2_V\}.
\end{align}
Moreover, similarly to \eqref{V-N}, there exist $\kappa, C \in
(2,\9)$ such that
\begin{align} \label{V'-N}
   \min \{ \|u\|_{V'},\|u\|^{\kappa}_{V'} \}
   \leq C \int \wt{N}(|u|) d\xi.
\end{align}
(See the Appendix for a proof.) Then, \eqref{wtN-V} and \eqref{V'-N}
imply that
\begin{align} \label{bdd-v'}
    \| I_{\{|X_n|\leq e^{-3}\}} X_n L(|X_n|^2) \|_{V'}
    \leq& C \max \{ \|X_n\|_{V},\|X_n\|^2_{V}, \|X_n\|^{1/\kappa}_V, \|X_n\|^{2/\kappa}_V  \} \nonumber \\
    \leq& C(\|X_n\|_V^2 + 1),
\end{align}
Hence, since $p\geq 3$, $2p'\leq p$, H\"{o}lder's inequality yields
\begin{align} \label{bdd-lp'-v'}
    \|I_{\{|X_n| \leq e^{-3}\}}  X_n L(|X_n|^2)   \|_{L^{p'}(\Omega
    \times(0,T);V')}
    \leq& C_T(\|X_n\|^2_{L^{2p'}(\Omega\times (0,T); V)}+1)
    \nonumber\\
    \leq& C_T (\|X_n\|^2_{L^{p}(\Omega\times (0,T); V)}+1) \nonumber\\
    \leq& C_T<\9,
\end{align}
where $C_T$ is independent of $n$.

Consequently, \eqref{split-u'}, \eqref{bdd-lp'-h-1} and
\eqref{bdd-lp'-v'} together yield \eqref{bdd-lp'-u'}, thereby
completing the proof of Proposition \ref{Max-ABG}. \hfill $\square$

\begin{remark}
As in \cite{BR14}, we can also define the operators $\calb, \cala:
\mathcal{U} \to \mathcal{U}'$ by
\begin{align*}
  (\cala y)(t)=& ie^{-W(t)}\D (e^{W(t)}y(t)) + 4i\lbb|\lbb|t y(t),\ y\in D(\cala)=\mathcal{U}, \nonumber  \\
  (\calb y)(t)=&\frac{dy(t)}{dt}+ \wh{\mu}y(t),\ a.e.\ t\in(0,T),\ y\in D(\calb),
\end{align*}
where $
  D(\calb) = \{y\in \mathcal{U}:\  y\in AC([0,T];U')\cap C([0,T];H),
  \bbp-a.s.,  \frac{dy}{dt} \in \mathcal{U}', y(0)=x
  \}.$
Then, \eqref{equa-y} can be reformulated as an operatorial equation
\begin{align*}
   \calb y + \cala y + \calg y =0.
\end{align*}
It is clear that $\cala $ is maximal monotone from $\calu$ to
$\calu'$. The same assertion holds also for $\calb$, by similar
arguments as in the proof of \cite[Lemma 4.1, Lemma 4.2]{BR14}.
Then, since $D(\cala)=D(\calg)=\calu$, we deduce from \cite[Theorem
2.6]{B10} that $\cala + \calb + \calg$ is also maximal monotone.
However, unlike in \cite{BR14}, we do not have the coercivity (see
\cite[(2.3)]{BR14}) in the Schr\"{o}dinger case, the proof of
\cite[Proposition 3.3]{BR14} is not applicable here. In order to
obtain existence of solutions to \eqref{equa-y*}, we shall introduce
and study an associated approximating equation in the next section.
\end{remark}

\section{Approximating equation} \label{App-Equa}

Consider the approximating equation,
\begin{align} \label{app-equa-y}
   & y(t) = x - \int_0^t \(i e^{-W(s)} \D(e^{W(s)}y(s)) +
   (4i\lbb|\lbb|t  + \wh{\mu} )y(s) + \calg_\ve (y(s))\)ds,\\
   &y(0)=x, \nonumber
\end{align}
Here, $t\in (0,T)$, $0\leq \ve \leq 1$,
\begin{align} \label{def-calg}
    \calg_{\ve} (y)
    :=2 \lbb i y  L_{\ve}(e^Wy)
     + 2|\lbb|y,
\end{align}
and
\begin{align} \label{def-l}
    L_{\ve}(u) = \log (\frac{|u| +\ve}{1+\ve |u|}),\ \forall u\in \bbc.
\end{align}
For $\ve =0$, set $L(u):= L_0(u) = \log |u|$, $u\in \bbc$.

We collect some properties of $L_\ve$ in the following lemma, whose
proof is included in the Appendix for completeness.

\begin{lemma} \label{Pro-Lve}
Let $0<\ve < 1$. Then:

$(i)$ For all $u>0$, $|L_\ve (u)| \leq |\log \ve|$, and $|u L_\ve
(u)| \leq |u L(u)|$.

$(ii)$ For all $u_1, u_2 \in \mathbb{C}$,
\begin{align} \label{Lve-diff}
    \left| u_1L_{\ve}(u_1)  - u_2L_{\ve}(u_2) \right|
    \leq& (1+\log(1/\ve))|u_1-u_2|.
\end{align}

$(iii)$ For all $u_1,u_2\in \mathbb{C}$,
\begin{align} \label{mono-lve-0}
   | Im (\ov{u_1} - \ov{u_2}) (u_1 L_\ve(u_1) - u_2 L_\ve(u_2)) |
   \leq (1-\ve^2 ) |u_1 - u_2|^2.
\end{align}
\end{lemma}

The main result in this section is as follows.

\begin{proposition} \label{Prop-yve-op}
Assume $(H)$ and let $0<\ve<1$ be fixed.  For any initial datum
$x\in U$ and $0<T<\9$, there exists a unique $U$-valued
$(\mathscr{F}_t)$-adapted  process $y_\ve$, such that $y_\ve\in
C([0,T];H^1)$, $\bbp$-a.s., and it satisfies \eqref{app-equa-y} in
the space $U'$ on $[0,T]$, $\bbp$-a.s.

Moreover, for any $p\geq 2$,
\begin{align} \label{esti-u-p}
   \bbe \sup\limits_{0\leq t\leq T} \|e^{W(t)}y_\ve(t)\|^p_U
   \leq C(T,p)<\9,
\end{align}
and
\begin{align} \label{esti-u'-p}
    \bbe \sup\limits_{0\leq t\leq T} \|e^{W(t)}
    \calg_\ve(y_\ve(t))\|^p_{U'}
   \leq C(T,p)<\9,
\end{align}
where $C(T,p)$ is independent of $\ve$.
\end{proposition}

The proof will proceed in two steps. We first prove the global
well-posedness of \eqref{app-equa-y} in the state space  $H^1$  in
Subsection \ref{H1-GWP}, and then we prove the necessary uniform
estimates in the Orlicz space in Subsection \ref{Uni-Log}.

\subsection{$H^1$ global well-posedness} \label{H1-GWP}

\begin{proposition} \label{Prop-y-h1}
Assume $(H)$ and let $0<\ve<1$ be fixed. For each $x\in H^1$ and
$0<T<\9$, there exists a unique $H^1$-valued
$(\mathscr{F}_t)$-adapted process $y_\ve$, such that $y_{\ve} \in
C([0,T]; H^1)$, and it solves \eqref{app-equa-y} in the space
$H^{-1}$ on $[0,T]$, $\bbp$-a.s.

Moreover, for any $p\geq 2$,
\begin{align} \label{esti-bdd-h1-gloabl}
   \bbe \sup\limits_{t\in[0,T]} |e^{W(t)}y_{\ve}(t)|_{H^1}^p \leq C(T,p) <\9,
\end{align}
where $C(T,p)$ is independent of $\ve$.
\end{proposition}

The key observation for the proof lies in the fact that the operator
$y\to y L_{\ve}(e^Wy)$ is Lipschitz on $L^2$ and bounded on $H^1$.
This fact allows to apply  a fixed point argument as in
\cite{BRZ14.2}. Below, the proof will rely on three lemmas. We first
introduce the evolution operators in Lemma \ref{Evolu-A}, and then
we prove the local existence in Lemma \ref{Local-ye}. Finally, in
Lemma \ref{Esti-h1-energy} we derive the a priori estimate in
$H^1$-norm, which in turn implies the global well-posedness.

\begin{lemma} \label{Evolu-A}
$\mathbb{P}-a.e.$, the operator $y\to -ie^{-W}\Delta (e^Wy) - (
2|\lbb| + 4i\lbb|\lbb|t + \wh{\mu}) y$ generates evolution operators
$U(t,s)=U(t,s,\omega)$ in the space $H^1(\mathbb{R}^d)$, $0\leq
s\leq t\leq T$. For each $x\in H^1(\mathbb{R}^d)$ and $s\in[0,T]$,
the process $[s,T]\ni t \to U(t,s)x$ is continuous and
$(\mathscr{F}_t)$-adapted, hence progressively measurable with
respect to the filtration $(\mathscr{F}_t)_{t\geq s}$.

Moreover, for any $f\in L^1(0,T; H^1)$, then $H^1$-path
\begin{align} \label{equa-U}
   y(t) =U(t,0)x + \int_0^t U(t,s) f(s) ds, \ \ 0\leq t\leq T,
\end{align}
satisfies the estimates
\begin{align} \label{Esti-U-0}
   \|y\|_{C([0,T]; H)} \leq C_T (|x|_{H} + \|f\|_{L^1(0,T;
   H)}),
\end{align}
and
\begin{align} \label{Esti-U}
   \|y\|_{C([0,T]; H^1)} \leq C_T (|x|_{H^1} + \|f\|_{L^1(0,T;
   H^1)}).
\end{align}
Here, the process $C_t$, $t\geq 0$, can be taken to be
$(\mathscr{F}_t)$-adapted progressively measurable, increasing and
continuous.
\end{lemma}
(For the proof see the Appendix.)

\begin{lemma} \label{Local-ye}
Assume $(H)$ and let $0<\ve<1$ be fixed. For each $x\in H^1$ and
$0<T<\9$, there exists an $H^1$-valued $(\mathscr{F}_t)$-adapted
process $y_{\ve}$ and a stopping time $\tau_{\ve}^*(x) \leq T$, such
that $y_\ve \in C([0,\tau_{\ve}^*(x)); H^1)$, and $y_\ve$ solves the
equation \eqref{app-equa-y} in $H^{-1}$  on $[0,\tau_{\ve}^*(x) )$,
$\bbp$-a.s.

Moreover, $\tau^*_\ve(x) =T$, $\bbp-a.s$, if
\begin{align} \label{yve-blow}
   \sup\limits_{t\in[0,\tau_\ve^*(x))} |y_\ve(t)|_{H^1}<\9,\ \  \bbp-a.s.
\end{align}
\end{lemma}

{\it Proof.}  Using the evolution operators introduced in Lemma
\ref{Evolu-A}, we reformulate the equation \eqref{app-equa-y} in the
mild form
\begin{align} \label{mild-equay-ve}
   y= U(t,0)x -2\lbb i \int_0^t U(t,s)\( y(s) L_{\ve} (e^{W(s)} y(s))\)ds,
\end{align}
(Note that, since for $y\in C([0,T];H^1)$, $yL_\ve(e^Wy)\in
L^1(0,T;H^1)$, the equivalence between \eqref{app-equa-y} and
\eqref{mild-equay-ve} can be proved similarly as in \cite[Theorem
2.2.2]{Z14}.)

Consider the integral operator $F$ defined for any $y \in C([0,T] ;
H^1)$ by
\begin{align*}
  F(y)(t) :=   U(t,0)x -2 \lbb i \int_0^t U(t,s) ( y(s) L_{\ve}
(e^{W(s)}
  y(s)))ds,\ t\in[0,T].
\end{align*}
We first show that
\begin{align} \label{F-h1-bdd}
   F(C([0,T] ; H^1)) \subset C([0,T] ; H^1).
\end{align}
Indeed, by \eqref{Esti-U},
\begin{align*}
  \|F(y)\|_{C([0,T] ; H^1)}
  \leq& C_T \(|x|_{H^1}
       + 2|\lbb| \|y L_{\ve}(e^W y)\|_{L^1(0,T;H^1)}
       \).
\end{align*}
By Lemma \ref{Pro-Lve} $(i)$ we have
\begin{align*}
   |yL_{\ve}(e^Wy)|_{H^1}
   \leq& \sqrt{2} |\log \ve|  |y|_{H^1} + |y \na (L_{\ve}(e^Wy))|_2.
\end{align*}
Moreover, straightforward computations show that
\begin{align} \label{cal-na-le}
   \na (L_{\ve}(e^Wy))
   = \frac{(1-\ve^2) |e^Wy|^{-1} Re(\ov{e^Wy} \na (e^Wy))}{(\ve+ |e^Wy|)(1+ \ve
   |e^Wy|)},
\end{align}
which implies that
\begin{align} \label{esti-nal}
   |\na (L_\ve(e^Wy))| \leq |e^Wy|^{-1} |\na(e^Wy)|.
\end{align}
Then,
\begin{align*}
   |y\na (L_{\ve}(e^Wy))|_2 \leq |e^{-W} \na (e^Wy)|_2
   \leq  \sqrt{2}\exp(2|W|_{L^\9})(1+ |\na W|_{L^\9}) |y|_{H^1}.
\end{align*}
Hence,
\begin{align*}
  |y L_{\ve}(e^Wy)|_{H^1}
  \leq \sqrt{2}( |\log \ve|  +  \exp(2|W|_{L^\9})(1+ |\na W|_{L^\9}) ) |y|_{H^1}.
\end{align*}
It follows that
\begin{align} \label{esti-F-bdd}
  \|F(y)\|_{C([0,T];H^1)}
  \leq C_T |x|_{H^1}
       + C_T D_{1}(T) T \|y\|_{C([0,T];H^1)}
\end{align}
with $D_{1}(T):=  2\sqrt{2}|\lbb|( |\log \ve|  + \sup\limits_{t\leq
T}\exp(2|W(t)|_{L^\9})(1+ \sup\limits_{t\leq T}|\na W(t)|_{L^\9}) )
$, thereby yielding \eqref{F-h1-bdd} as claimed.

Next, we will apply the iteration arguments as in \cite{BRZ14.2} to
construct the local solution to \eqref{app-equa-y}.

Fix $\oo \in \Omega$. Set $\mathcal {Y}^{\tau_1}_{M_1}:= \{y\in
C([0,\tau_1]; H^1): \|y\|_{C([0,\tau]; H^1)} \leq M_1\}$, where
$\tau_1$ and $M_1$ are random variables to be chosen later.

Similarly to \eqref{esti-F-bdd}, for any $y\in \mathcal
{Y}^{\tau_1}_{M_1}$,
\begin{align} \label{esti-F-bdd-tau}
  \|F(y)\|_{C([0,\tau_1];H^1)}
  \leq C_{\tau_1} |x|_{H^1}
       + C_{\tau_1} D_{1}(\tau_1) M_1  \tau_1.
\end{align}
Moreover,  for any $y, \wt{y} \in \mathcal {Y}^{\tau_1}_{M_1}$, by
\eqref{Esti-U-0},
\begin{align*}
   \|F(y) - F(\wt{y})\|_{C([0,\tau_1];L^2)}
   \leq 2|\lbb| C_{\tau_1}
         \|yL_{\ve}(e^Wy) - \wt{y}L_{\ve}(e^W\wt{y})\|_{L^1(0,\tau_1;
          L^2)},
\end{align*}
which implies by \eqref{Lve-diff} that
\begin{align}  \label{esti-F-lip-tau}
   \|F(y) - F(\wt{y})\|_{C([0,\tau_1];L^2)}
   \leq C_{\tau_1} D_{2}(\tau_1) \tau_1
   \|y-\wt{y}\|_{C([0,\tau_1];L^2)},
\end{align}
where $D_{2}(t)= 2|\lbb| (1+ |\log \ve|  ) \sup\limits_{s\leq
t}\exp(2|W(s)|_{L^\9})$.

Then, we define the real-valued continuous,
$(\mathscr{F}_t)$-adapted process $Z(t):=D_1(t)+D_2(t)$, and denote
the $(\mathscr{F}_t)$-stopping time $\tau_1 := \inf\{t\in[0,T]: C_t
Z(t)t\geq \frac 12 \}\wedge T $ and $M:= 2 C_{\tau_1} |x|_{H^1}$.
\eqref{esti-F-bdd-tau} and \eqref{esti-F-lip-tau} imply that
$F(\mathcal{Y}^{\tau_1}_{M_1}) \subset \mathcal{Y}^{\tau_1}_{M_1}$
and $F$ is a contraction in $C([0,\tau_1];L^2)$. Hence, Banach's
fixed point theorem yields a unique $y\in
\mathcal{Y}^{\tau_1}_{M_1}$, such that $y = F(y)$ on $[0,\tau_1]$.
Setting $y_1(t):= y(t\wedge \tau_1)$ and arguing as in the proof of
\cite[Proposition $2.5$]{BRZ14.2} we deduce that
$y_{1}|_{[0,\tau_1]}\in C([0,\tau_1]; H^1)$, $y_1$ is
$(\mathscr{F}_t)$-adapted and it solves \eqref{app-equa-y} on
$[0,\tau_1]$, $\bbp$-a.s.

Applying similar arguments as in \cite{BRZ14.2}, we can extend the
solution step by step and construct a sequence $\{(y_m,
\tau_m)\}_{m\geq 1}$, such that for each $m\geq 1$, $\tau_m$ is an
$(\mathscr{F}_t)$-stopping time, $\tau_{m+1} \geq \tau_m$, $y_m$ is
an $H^1$-valued $(\mathscr{F}_t)$-adapted process, such that
$y_m|_{[0,\tau_m]} \in C([0,\tau_m]; H^1)$, $y_m(t) = y_m(t \wedge
\tau_m)$, $t\in [0,T]$, and $y_m$ solves \eqref{app-equa-y} on
$[0,\tau_m]$, $\bbp$-a.s.

More precisely, given the pair $(y_m,\tau_m)$ with such properties
above at the $m$-th step, we set $\mathcal
{Y}^{\sigma_m}_{M_{m+1}}:= \{ z\in C([0,\sigma_m]; H^1):
\|z\|_{C([0,\sigma_m]; H^1)} \leq M_{m+1}\}$, and define for $z\in
C([0,T];H^1)$,
\begin{align*}
    F_m(z)(t) := &  U(\tau_m+t,\tau_m)y_m(\tau_m) \\
       & - 2\lbb i \int_0^t U(\tau_m+t,\tau_m+s) ( z(s) L_{\ve}
(e^{W(\tau_m+s)} z(s)))ds.
\end{align*}
Similarly to \eqref{esti-F-bdd-tau} and \eqref{esti-F-lip-tau}, for
$z\in \mathcal {Y}^{\sigma_m}_{M_{m+1}}$,
\begin{align*}
    \|F_m(z)\|_{C([0,\sigma_m];H^1)}
  \leq C_{\tau_m+\sigma_m} (|y_m(\tau_m)|_{H^1}
       +  D_{1}(\tau_m+\sigma_m) M_{m+1}\sigma_m).
\end{align*}
and for $z,\wt{z} \in \mathcal {Y}^{\sigma_m}_{M_{m+1}}$,
\begin{align*}
    \|F_m(z) - F_m(\wt{z})\|_{C([0,\sigma_m];L^2)}
   \leq C_{\tau_m+\sigma_m} D_2(\tau_m+\sigma_m) \sigma_m
   \|z-\wt{z}\|_{C([0,\sigma_m];L^2)}.
\end{align*}
Then, define $Z_t^{(m)}:= D_1(\tau_m+t) + D_2(\tau_m+t)$ and
$\sigma_m = \inf \{t\in[0,T-\tau_m]: C_{\tau_m+t} Z_t^{(m)} t>\frac
12 \} \wedge (T-\tau_m)$. It follows that $F_m(\mathcal
{Y}^{\sigma_m}_{M_{m+1}}) \subset \mathcal {Y}^{\sigma_m}_{M_{m+1}}$
and $F_m$ is a contraction in $C([0,\sigma_m];L^2)$. By Banach's
fixed point theorem, we obtain a unique $z_{m+1}\in \mathcal
{Y}^{\sigma_m}_{M_{m+1}}$, such that $z_{m+1} = F_m(z_{m+1})$ on
$[0,\sigma_m]$.

Therefore, set $\tau_{m+1}: = \tau_m + \sigma_m$ and
\begin{align*}
    y_{m+1}(t) :=
    \left\{
      \begin{array}{ll}
        y_m(t), & \hbox{$t\in[0,\tau_m]$;} \\
        z_{m+1}((t-\tau_m)\wedge \sigma_m), & \hbox{$t\in(\tau_m,T]$.}
      \end{array}
    \right.
\end{align*}
Then, we construct a new pair $(y_{m+1},\tau_{m+1})$ with the
properties mentioned above. In particular, $y_{m+1}$ solves
\eqref{app-equa-y} on $[0,\tau_{m+1}]$, $\bbp$-a.s. Iterating this
procedure gives us the desired sequence $\{(y_m,\tau_m)\}_{m\geq
1}$.

Now, let $\tau_\ve^*(x): =\lim\limits_{m\to \9} \tau_m$ and
$y_{\ve}:= \lim\limits_{m\to \9} y_m I_{[0,\tau_\ve^*(x))}$. It
follows that $\tau_\ve^*(x)$ is an $(\mathscr{F}_t)$-stopping time,
$y_\ve$ is an $H^1$-valued $(\mathscr{F}_t)$-adapted process,
$y_\ve\in C([0,\tau_\ve^*(x)); H^1)$, and it solves the equation
\eqref{app-equa-y} on $[0,\tau_\ve^*(x))$, $\bbp$-a.s.

Finally, by the construction of $\{(y_m,\tau_m)\}_{m\geq 1}$, we use
similar arguments as in \cite{BRZ14.2} to obtain the blow-up
alternative, i.e. for $\bbp$-a.e. $\omega$, if $\tau_m(\omega) <
\tau^*_\ve(x)(\omega)$, $\forall m\in \mathbb{N}$, then $
    \lim\limits_{t\to \tau_\ve^*(x)(\omega)} |y_\ve(t)(\omega)|_{H^1}
=\9.$ By the construction of $\sigma_m$ above, we consequently
conclude that $\tau_\ve^*(x)=T$ if \eqref{yve-blow} holds. \hfill
$\square$

\begin{lemma} \label{Esti-h1-energy}
Assume the conditions of Lemma \ref{Local-ye} to hold, and let
$\tau_\ve^*(x)$ and $y_\ve$ be as in Lemma \ref{Local-ye}. Then, for
any $p\geq 2$,
\begin{align} \label{esti-bdd-h1-localye}
   \bbe \sup \limits_{t\in[0,\tau_\ve^*(x))}
|e^{W(t)}y_{\ve}(t)|^p_{H^1} \leq C(T,p) < \9,
\end{align}
where $C(T,p)$ is independent of $\ve$.
\end{lemma}

{\it Proof.} Let $X_m:=e^Wy_\ve$, $\phi_j=\mu_je_j$, $1\leq j\leq
m$, and $\{\tau_m\}_{m\geq 1}$ be the sequence of stopping times
constructed in the proof of Lemma \ref{Local-ye}. Since $X_\ve
L_\ve(X_\ve)\in L^2 \subset H^{-1}$, as in the proof of \cite[Lemma
$5.2$]{BRZ14.2}, we derive that $\bbp$-a.s., for $t\in[0,\tau_m]$,
\begin{align} \label{ito-app-equa-xa}
    & | X_{\ve}(t)|_{H^1}^2 \nonumber \\
   =& |x|_{_{H^1}}^2
       -4|\lbb| \int_0^t |X_\ve|_{H^1}^2 ds
     -2 \int_0^t Re \int \na  \ov{X_{\ve}} \na (\mu X_{\ve}) d\xi ds
       \nonumber \\
     &  + \sum\limits_{j=1}^n \int_0^t |\na (X_{\ve}\phi_j)|_{2}^2
    ds   + 4\lbb \int_0^t Im \int \ \na \ov{ X_{\ve} }\na (X_{\ve}L_{\ve}(X_{\ve})) d\xi
    ds\nonumber \\
     &+ 2 \sum\limits_{j=1}^n \int_0^t \int |X_\ve|^2 Re \phi_j d\xi d\beta_j(s)
     +2 \sum\limits_{j=1}^n \int_0^t Re \int \na \ov{ X_{\ve}}
     \na (X_{\ve}\phi_j) d\xi d\beta_j(s).
\end{align}
Then, applying It\^{o}'s formula we obtain for any $p\geq 2$,
\begin{align}
  &|X_\ve(t)|_{H^1}^p \nonumber \\
  =& |x|_{H^1}^p -2p|\lbb| \int_0^t |X_\ve|_{H^1}^p ds \nonumber \\
  & -p\int_0^t |X_\ve|_{H^1}^{p-2}  Re\int \na\ov{X_\ve} \na(\mu
  X_\ve) d\xi ds
   + \frac{p}{2} \sum\limits_{j=1}^n \int_0^t |X_\ve|_{H^1}^{p-2} |\na
  (X_\ve\phi_j)|_2^2 ds \nonumber \\
  & + \frac{1}{2}p(p-2) \sum\limits_{j=1}^n \int_0^t |X_\ve|_{H^1}^{p-4}
      \(\int|X_\ve|^2 Re \phi_j d\xi + Re \int \na\ov{X_\ve}
      \na(X_\ve\phi_j)d\xi \)^2 ds \nonumber\\
  &+ 2p\lbb \int_0^t |X_\ve|_{H^1}^{p-2} Im \int \na\ov{X_\ve} \na(X_\ve
  L_\ve(X_\ve))d\xi ds\nonumber \\
  &+ p \sum\limits_{j=1}^n \int_0^t |X_\ve|_{H^1}^{p-2} \int |X_\ve|^2 Re \phi_j d\xi
  d\beta_j(s)\nonumber \\
  &+ p \sum\limits_{j=1}^n \int_0^t |X_\ve|_{H^1}^{p-2} Re \int \na\ov{X_\ve} \na
  (X_\ve\phi_j) d\xi d\beta_j(s)\nonumber
\end{align}
\begin{align}
  =:&|x|_{H^1}^p + \sum\limits_{k=1}^7 J_k(t),\ \ t\in[0,\tau_m],\
  \bbp-a.s.
\end{align}

Since $e_j\in C_b^{\9}$, $1\leq j\leq n$ and since by
\eqref{esti-nal} we have
\begin{align*}
   \bigg| Im \int  \na \ov{X_{\ve}} \na ( X_{\ve}L_{\ve}(X_{\ve})) d\xi \bigg|
   = \bigg| Im\ \int \na  \ov{ X_\ve} (X_\ve \na L_\ve(X_\ve)) d\xi
\bigg| \leq |\na X_\ve|_2^2,
\end{align*}
it follows that
\begin{align} \label{esti-J-15}
   \sum\limits_{k=1}^5 \bbe \sup\limits_{s\leq t} |J_k(s)|
   \leq C_p \int_0^t \bbe \sup\limits_{r\leq s}
   |X_\ve(r)|^p_{H^1} ds,\ \ t\in[0,\tau_m],
\end{align}
where $C_p$ is independent of $\ve$ and $m$.

As regards the remaining stochastic terms, it follows from the
Burkholder-Davis-Gundy inequality that
\begin{align} \label{esti-J-6}
   \bbe \sup\limits_{s\leq t\wedge \tau_m} |J_6(s)|
   \leq& C_p \bbe \bigg[ \int_0^{t\wedge \tau_m} \sum\limits_{j=1}^n |X_\ve|^{2p-4}_{H^1} \(\int |X_\ve|^2 Re \phi_j d\xi\)^2 ds \bigg]^{\frac
   12} \nonumber \\
   \leq& C_p \bbe \( \int_0^{t\wedge \tau_m} |X_\ve|^{2p}_{H^1} ds\)^{\frac 12}
   \nonumber \\
   \leq& C_p \delta \bbe \sup\limits_{s\leq t\wedge \tau_m} |X_\ve|^p_{H^1}
         + C(p,\delta) \int_0^t \bbe \sup\limits_{r\leq s\wedge \tau_m}
         |X_\ve|^p_{H^1} ds,
\end{align}
where we used \cite[Lemma $3.3$]{BRZ14.2} in the last step,
$\delta>0$, and $C_p,C(p,\delta)$ are independent of $\ve$ and $m$.

Similarly,
\begin{align} \label{esti-J-7}
   \bbe \sup\limits_{s\leq t\wedge \tau_m} |J_7(s)|
   \leq& C_p \bbe \bigg[ \int_0^{t\wedge \tau_m} \sum\limits_{j=1}^n  |X_\ve|^{2p-4}_{H^1} \(Re \int \na\ov{X_\ve} \na
   (X_\ve\phi_j)d\xi\)^2 ds \bigg]^{\frac 12} \nonumber \\
   \leq& C_p\delta \bbe \sup\limits_{s\leq t\wedge \tau_m} |X_\ve|^p_{H^1}
         + C(p,\delta) \int_0^t \bbe \sup\limits_{r\leq s\wedge \tau_m}
         |X_\ve|^p_{H^1} ds,
\end{align}
where $C(p,\delta)$ is independent of $\ve$ and $m$.

Therefore, combining \eqref{esti-J-15}-\eqref{esti-J-7}, taking
$\delta$ sufficiently small, and applying the Gronwall inequality we
obtain
\begin{align*}
   \bbe \sup\limits_{t\in[0,\tau_m]} |X_\ve(t)|^p_{H^1} \leq
   C(T,p)<\9,
\end{align*}
where $C(T,p)$ is independent of $\ve$ and $m$. Taking $m\to \9$ and
using Fatou's lemma we consequently
obtain \eqref{esti-bdd-h1-localye}. \hfill $\square$\\

{\it \bf Proof of Proposition \ref{Prop-y-h1}.} It follows from
\eqref{esti-bdd-h1-localye}  that
$\sup\limits_{t\in[0,\tau^*_\ve(x))} |e^{W(t)}y_\ve(t)|_{H^1} <\9$,
$\bbp$-a.s.  Then, since $e_j\in C_b^\9$, $1\leq j\leq n$,
$\sup\limits_{t\in[0,\tau^*_\ve(x))} |y_\ve(t)|_{H^1} <\9$,
$\bbp$-a.s., which along with Lemma \ref{Local-ye} implies the
global existence of the solution to \eqref{app-equa-y}.

Uniqueness for  \eqref{app-equa-y} follows from monotonicity.
Indeed, consider any two solutions $y_1, y_2$ to \eqref{app-equa-y}
with the initial datum $x$, and set $X_i =e^W y_i$, $i=1,2$. Then,
similarly to \eqref{ito-app-equa-xa}, we derive that
\begin{align} \label{difference}
   \bbe |X_1(t) - X_2(t)|_2^2
   =& -4   |\lbb|  \bbe \int_0^t |X_1 - X_2|_2^2 ds \nonumber \\
    &+ 4\lbb  \bbe  \int_0^t Im  \int (\ov{X_1}- \ov{X_2})(X_1 L_\ve(X_1) - X_2 L_\ve(X_2)
    ) d\xi ds.
\end{align}
By \eqref{mono-lve-0},
\begin{align} \label{mono-lve}
    &\bigg| \bbe  \int_0^t Im \int  (\ov{X_1}- \ov{X_2})(X_1 L_\ve(X_1) - X_2 L_\ve(X_2)
    ) d\xi ds \bigg| \nonumber \\
    \leq& (1-\ve^2)\bbe \int_0^t |X_1-X_2|_{2}^2 ds \leq  \bbe \int_0^t|X_1-X_2|_2^2ds.
\end{align}
Then, it follows that
\begin{align*}
   \bbe |X_1(t) - X_2(t)|_2^2 \leq 0,
\end{align*}
which implies that for each $t\in[0,T]$,  $X_1(t)= X_2(t)$,
$\bbp$-a.s. Thus, by the continuity of $y_i$ in $H^1$, $i=1,2$, we
deduce that $X_1(t)= X_2(t)$, $\forall t \in [0,T]$, $\bbp$-a.s.,
thereby obtaining the uniqueness. \hfill $\square$ \\

In the next subsection, we shall derive some uniform estimates in
the Orlicz space, which allows to apply the method of maximal
monotone operators to take the limit in the approximating equations.

\subsection{Uniform estimates} \label{Uni-Log}

This subsection is mainly devoted to uniform estimates in the Orlicz
space $V$. Taking into account the definition \eqref{def-V}, let us
begin with the estimate of the entropy function below.

\begin{lemma} \label{Esti-Log}
Fix $0<\ve \leq 1$. Let $x \in U$, $0<T<\9$, $y_\ve$ be the
approximating solution in Proposition \ref{Prop-y-h1} and
$X_\ve=e^Wy_\ve$. We have for any $p\geq 2$,
\begin{align} \label{esti-log}
  \bbe \sup\limits_{t\leq T} \bigg| \int |X_{\ve}(t)|^2 \log |X_{\ve}(t)|^2
  d\xi \bigg|^p \leq C(T,p)<\9,
\end{align}
where $C(T,p)$ is independent of $\ve$.
\end{lemma}

{\it \bf Proof.} For $u>0$, set $F_m(u):= \int_0^u (L_{1/m}(\nu) +
1) d\nu$, where $L_{1/m}(\cdot)$ is as defined in \eqref{def-l}.
Using the techniques as in \cite{K10} and \cite[Lemma 5.1]{BRZ14.2}
we can derive that $\bbp$-a.s., for $t\in [0,T]$,
\begin{align} \label{Ito-Fn}
    &\int F_m(|X_\ve(t)|^2) d\xi  \nonumber \\
    =& \int F_m(|x|^2) d\xi
       - 2\int_0^t \int g_m(|X_{\ve}|^2) Im(\ov{X_{\ve}} \na
       X_{\ve})
       Re (\ov{X_{\ve}} \na X_{\ve}) d\xi ds  \nonumber \\
     & -4|\lbb| \int_0^t \int (L_{1/m}(|X_\ve|^2)+1)|X_\ve|^2 d\xi ds
       + \int_0^t \int  g_m(|X_{\ve}|^2) (Re \phi_j)^2|X_{\ve}(s)|^4 d\xi ds  \nonumber \\
     & + 2 \sum\limits_{j=1}^n \int_0^t \int (L_{1/m}(|X_\ve|^2)+1)|X_{\ve}|^2 Re \phi_j d\xi
       d\beta_j(s),
\end{align}
where  $g_m (|X_{\ve}|^2):= 2(1-m^{-2}) (m^{-1}
+|X_{\ve}|^2)^{-1}(1+ m^{-1}|X_{\ve}|^2)^{-1}$, and $\phi_j=\mu_j
e_j$, $1\leq j\leq n$. (See the Appendix for the proof.)

Then, applying It\^{o}'s formula we derive that $\bbp$-a.s. for
$t\in [0,T]$,
\begin{align} \label{Ito-Fn-p}
   &\(\int F_m(|X_\ve(t)|^2) d\xi \)^p \nonumber \\
   =& \(\int F_m(|x|^2) d\xi \)^p \nonumber \\
    &-2 p \int_0^t \(\int F_m(|X_\ve|^2) d\xi \)^{p-1}
       \(\int g_m(|X_\ve|^2) Im(\ov{X_\ve}\na X_\ve) Re  (\ov{X_\ve}\na
       X_\ve) d\xi\) ds \nonumber \\
    & -4 |\lbb| p \int_0^t \(\int F_m(|X_\ve|^2) d\xi \)^{p-1}
     \(\int\(L_{1/m}(|X_\ve|^2)+1\)|X_\ve|^2 d\xi \) \ ds \nonumber \\
   &+ p \int_0^t \(\int F_m(|X_\ve|^2) d\xi \)^{p-1}
        \(\int g_m(|X_\ve|^2) (Re \phi_j)^2 |X_\ve|^4 d\xi\) ds
        \nonumber \\
    & +2 p(p-1) \sum\limits_{j=1}^N \int_0^t \(\int F_m(|X_\ve|^2) d\xi \)^{p-2}
         \bigg(\int (L_{1/m}(|X_\ve|^2)+1)|X_\ve|^2 Re \phi_j d\xi\bigg)^2ds
         \nonumber
\end{align}
\begin{align}
    & + 2p \sum\limits_{j=1}^n\int_0^t \(\int F_m(|X_\ve|^2)d\xi \)^{p-1} \(\int |X_\ve|^2 Re
    \phi_j d\xi\) d\beta_j(s) \nonumber \\
    & + 2p \sum\limits_{j=1}^n \int_0^t \(\int F_m(|X_\ve|^2)d\xi \)^{p-1}
           \(\int |X_\ve|^2 L_{1/m}(|X_\ve|^2)  Re\phi_j d\xi\)
           d\beta_j(s). \nonumber \\
    =& \(\int F_m(|x|^2) d\xi \)^p + \sum\limits_{j=1}^6 K_j(t).
\end{align}

Since for $u>0$,
\begin{align} \label{Fm-L}
   F_m(u) = uL_{1/m}(u)+ u - (1-m^{-2})\int_0^u
   \nu(m^{-1}+\nu)^{-1}(1+m^{-1}\nu)^{-1} d\nu ,
\end{align}
it follows that
\begin{align*}
   |F_m(|x|^2)| \leq ||x|^2 L_{1/m}(|x|^2)| + 2|x|^2
   \leq ||x|^2 L(|x|^2)| + 2|x|^2 \in L^1(\bbr^d).
\end{align*}
Then, since $F_m(|x|^2) \to |x|^2 L(|x|^2)$, as $m\to \9$, the
dominated convergence theorem yields
\begin{align} \label{conv-Fm}
   \int F_m(|x|^2) d\xi \to \int |x|^2 L(|x|^2) d\xi.
\end{align}
In particular,
\begin{align} \label{bdd-Fm}
    \sup\limits_{m\geq 1} \bigg|\int F_m(|x|^2) d\xi \bigg| \leq C< \9.
\end{align}

For the other deterministic terms in \eqref{Ito-Fn-p}, since $e_j\in
C_b^\9$, $1\leq j\leq n$, and $g_m(|X_\ve|^2) \leq 2|X_\ve|^{-2}$,
using the Young inequality $a^{p-1}b\leq \frac{p-1}{p}\delta a^p +
\frac{1}{p} \delta^{-(p-1)}b^p$, $a^{p-2}b \leq \frac{p-2}{p} \delta
a^p + \frac{2}{p}\delta^{-\frac{p-2}{2}}b^{\frac p2}$, $\delta>0$,
and the boundedness of $H^1$-norm in \eqref{esti-bdd-h1-gloabl}, we
derive that $\bbp$-a.s.
\begin{align} \label{esti-K-15}
   &\sum\limits_{j=1}^4 \bbe \sup\limits_{s\leq t} |K_j(s)|
   \nonumber \\
   \leq& C_p T \delta \bbe \sup\limits_{s\leq t} \bigg|\int F_m(|X_\ve|^2) d\xi
   \bigg|^p \nonumber \\
    &+C(p,\delta) \bbe \sup\limits_{s\leq t} \bigg[ \int_0^s |X_\ve|^{2p}_{H^1} ds
         + \int_0^s \(\int \bigg||X_\ve|^2 L_{1/m}(|X_\ve|^2)\bigg|d\xi \)^p
dr  \bigg] \nonumber \\
   \leq& C(T,p,\delta)
      + C_p T \delta \bbe \sup\limits_{s\leq t} \bigg|\int F_m(|X_\ve|^2) d\xi
   \bigg|^p \nonumber \\
     &  + C(p,\delta) \bbe \sup\limits_{s\leq t}  \int_0^s \(\int \bigg||X_\ve|^2 L_{1/m}(|X_\ve|^2)\bigg|d\xi \)^p
   dr, \ t\in[0,T],
\end{align}
where $C_p, C(p,\delta)$ and $C(T,p,\delta)$ are independent of
$\ve$ and $m$.

Moreover, it follows from the Burkholder-Davis-Gundy inequality, the
Young inequality $a^{2p-2}b^2 \leq \frac{p-1}{p} a^{2p} +
\frac{1}{p}b^{2p}$, \eqref{esti-bdd-h1-gloabl} and Lemma $3.3$ in
\cite{BRZ14.2} that
\begin{align} \label{esti-K-6}
   &\bbe \sup\limits_{s\leq t} |K_5(s)| \nonumber \\
   \leq & C_p \bbe \bigg[ \int_0^t \sum\limits_{j=1}^n \bigg|\int F_m(|X_\ve|^2) d\xi
   \bigg|^{2p-2} \bigg| \int|X_\ve|^2 Re\phi_j d\xi \bigg|^2 ds
   \bigg]^{\frac 12} \nonumber \\
   \leq& C_p \bbe \bigg[ \int_0^t  \bigg|\int F_m(|X_\ve|^2) d\xi
   \bigg|^{2p}\bigg]^{\frac 12}
   +C_p \bbe \bigg[  \int_0^t \(\sum\limits_{j=1}^n  \bigg| \int|X_\ve|^2 Re\phi_j d\xi \bigg|^2\)^{p} ds
   \bigg]^{\frac 12} \nonumber \\
   \leq& C(T,p,\delta)+ C_p\delta \bbe \sup\limits_{s\leq t} \bigg|\int F_m(|X_\ve|^2) d\xi
   \bigg|^p
   + C(p,\delta) \int_0^t \bbe \sup\limits_{r\leq s} \bigg|\int F_m(|X_\ve|^2) d\xi
   \bigg|^pds.
\end{align}

Similarly,
\begin{align} \label{esti-K-7}
   &\bbe \sup\limits_{s\leq t} |K_6(s)|  \nonumber \\
   \leq& C_p \bbe \bigg[ \int_0^t  \bigg|\int F_m(|X_\ve|^2) d\xi
   \bigg|^{2p}\bigg]^{\frac 12} \nonumber \\
   & + C_p \bbe \bigg[  \int_0^t \(\sum\limits_{j=1}^n\bigg| \int |X_\ve|^2 L_{1/m}(|X_\ve|^2) Re \phi_j d\xi \bigg|^2\)^{p} ds   \bigg]^{\frac
   12} \nonumber \\
   \leq& C_p \delta \bbe \sup\limits_{s\leq t} \bigg[ \bigg|\int F_m(|X_\ve|^2) d\xi
   \bigg|^{p} + \(\int \bigg||X_\ve|^2 L_{1/m}(|X_\ve|^2) \bigg|d\xi  \)^p
   \bigg] \nonumber \\
   &+ C(p,\delta) \int_0^t \bbe \sup\limits_{r\leq s} \bigg[ \bigg|\int F_m(|X_\ve|^2) d\xi
   \bigg|^{p} + \(\int \bigg||X_\ve|^2 L_{1/m}(|X_\ve|^2) \bigg|d\xi  \)^p
   \bigg] ds
\end{align}

Thus, it follows from \eqref{bdd-Fm}-\eqref{esti-K-7} that
\begin{align} \label{Fm-delta}
    &\bbe \sup\limits_{s\leq t} \bigg| \int F_m(|X_\ve(s)|^2) d\xi
    \bigg|^p  \nonumber \\
    \leq& C(T ,p,\delta)
     + C(T,p)\delta \bbe \sup\limits_{s\leq t} \bigg[ \bigg| \int F_m(|X_\ve(s)|^2) d\xi
    \bigg|^p
     +  \(\int \bigg||X_{\ve}|^2 L_{1/m}(|X_{\ve}|^2) \bigg|d
    \xi\)^p \bigg]  \nonumber \\
    &+ C(T,p,\delta) \int_0^t  \bbe \sup\limits_{r\leq s}
        \bigg[\bigg| \int F_m(|X_\ve|^2) d\xi \bigg|^p
     +  \(
    \int
    \bigg| |X_{\ve}|^2 L_{1/m}(|X_\ve|^2)  \bigg| d\xi
    \)^p  \bigg] ds.
\end{align}
Since by \eqref{Fm-L},
\begin{align} \label{Fm-Lm}
   \bigg|\ \big| \int F_m(|X_\ve|^2) d\xi \big|
   - \big| \int |X_{\ve}|^2 L_{1/m}(|X_\ve|^2) d\xi \big|\ \bigg|
   \leq 2 |X_\ve|_2^2,
\end{align}
using \eqref{esti-bdd-h1-gloabl} we obtain
\begin{align} \label{esti-lm}
    &\bbe \sup\limits_{s\leq t} \bigg| \int |X_\ve|^2 L_{1/m}(|X_\ve|^2) d\xi
    \bigg|^p \nonumber \\
    \leq& C(T,p, \delta) + C(T,p)\delta \bbe \sup\limits_{s\leq t}\( \int \bigg||X_{\ve}|^2 L_{1/m}(|X_{\ve}|^2) \bigg| d
    \xi\)^p  \nonumber \\
    &+ C(T,p,\delta) \int_0^t  \bbe \sup\limits_{r\leq s} \(
    \int
    \bigg| |X_{\ve}|^2 L_{1/m}(|X_\ve|^2)  \bigg| d\xi
    \)^p ds.
\end{align}

Note that  $u^2 L_{1/m}(u^2) \leq u^2 \log u^2 \leq C_\delta(u^2 +
u^{2+\delta}) $ for $u>1$. The Sobolev imbedding theorem implies
that for $0<\delta<\frac{4}{d-2}$,
\begin{align} \label{log-x>1}
    \int I_{\{|X_\ve| > 1\}} |X_\ve|^2 L_{1/m}(|X_\ve|^2) d\xi
    \leq& C_\delta(|X_\ve|_{H^1}^{2} + |X_\ve|_{H^1}^{2+\delta}).
\end{align}
Hence,
\begin{align*}
     &\bigg| \int I_{\{|X_\ve|\leq 1\}} |X_\ve|^2 L_{1/m}(|X_\ve|^2) d\xi
     \bigg|^p \\
     =& \bigg|  \int|X_\ve|^2 L_{1/m}(|X_\ve|^2) d\xi
             - \int I_{\{|X_\ve| > 1\}} |X_\ve|^2 L_{1/m}(|X_\ve|^2) d\xi
             \bigg|^p \\
     \leq& C_p \bigg|  \int|X_\ve|^2 L_{1/m}(|X_\ve|^{2}) d\xi \bigg|^p
             + C(p,\delta) (|X_\ve|_{H^1}^{2p} +
             |X_\ve|_{H^1}^{(2+\delta)p}).
\end{align*}
and
\begin{align*}
     &\( \int \bigg| |X_\ve|^2 L_{1/m}(|X_\ve|^2)\bigg| d\xi  \)^p \\
     =& \(-\int  I_{\{|X_\ve|\leq 1\}} |X_\ve|^2 L_{1/m}(|X_\ve|^2) d\xi
              + \int I_{\{|X_\ve| > 1\}}   |X_\ve|^2 L_{1/m}(|X_\ve|^2)  d\xi
\)^p \\
     \leq& C_p \bigg|  \int I_{\{|X_\ve|\leq 1\}} |X_\ve|^2 L_{1/m}(|X_\ve|^2) d\xi \bigg|^p
             + C(p,\delta) (|X_\ve|_{H^1}^{2p} +
             |X_\ve|_{H^1}^{(2+\delta)p}).
\end{align*}

Therefore, inserting the two estimates above into \eqref{esti-lm}
and then using \eqref{esti-bdd-h1-gloabl} we get
\begin{align*}
   &\bbe \sup \limits_{s\leq t} \bigg| \int I_{\{|X_\ve|\leq 1\}}|X_\ve|^2 L_{1/m}(|X_\ve|^2) d\xi
   \bigg|^p \\
   \leq& C(T,p,\delta)
      + C(T,p) \delta \bbe \sup \limits_{s\leq t} \bigg| \int I_{\{|X_\ve|\leq 1\}}|X_\ve|^2 L_{1/m}(|X_\ve|^2) d\xi
   \bigg|^p \\
     & +  C(T,p,\delta) \int_0^t \bbe \sup \limits_{r\leq s} \bigg|\int I_{\{|X_\ve|\leq 1\}} |X_\ve|^2 L_{1/m}(|X_\ve|^2)
   d\xi \bigg|^p ds.
\end{align*}
Then, taking $\delta$ sufficiently small and applying Gronwall's
inequality we have
\begin{align*}
    \bbe \sup\limits_{t\leq T} \bigg| \int I_{\{|X_\ve|\leq 1\}} |X_{\ve}(t)|^2 L_{1/m}(|X_{\ve}(t)|^2) d\xi
    \bigg|^p
    \leq  C(T,p),
\end{align*}
where $C(T,p)$ is independent of $\ve$ and $m$. Hence, by Fatou's
lemma,
\begin{align} \label{esti-log-x1}
    \bbe \sup\limits_{t\leq T} \bigg| \int I_{\{|X_\ve|\leq 1\}} |X_{\ve}(t)|^2 L(|X_{\ve}(t)|^2) d\xi
    \bigg|^p
    \leq  C(T,p),
\end{align}

Consequently, \eqref{esti-log} follows immediately from
\eqref{esti-log-x1}, \eqref{log-x>1} and \eqref{esti-bdd-h1-gloabl}.
Hence, the proof is
complete. \hfill $\square$\\

{\it \bf Proof of Proposition \ref{Prop-yve-op}.} By Proposition
\ref{Prop-y-h1}, there exists a unique $(\mathscr{F}_t)$-adapted
solution $y_\ve$ to \eqref{app-equa-y}, and $y_\ve \in C([0,T];
H^1)$, $\bbp$-a.s. Moreover, since $u \mapsto u L_\ve(e^W u)$ is
Lipschitz on $L^2$ and $L^2 \subset U'$, we obtain $\calg_\ve
(y_\ve) \in C([0,T];U')$, $\bbp$-a.s.

It remains to prove \eqref{esti-u-p} and \eqref{esti-u'-p}. For the
proof of \eqref{esti-u-p}, in view of \eqref{esti-bdd-h1-gloabl}, we
only need to prove that for any $p\geq 2$,
\begin{align} \label{esti-v-p}
   \bbe \sup\limits_{0\leq t\leq T} \|X_\ve(t)\|^p_V
   \leq C(T,p)<\9,
\end{align}
where $X_\ve := e^W y_\ve$, and $V$ is the Orlicz space defined in
\eqref{def-V}.

To this end, set $B(u):= -u^2 \log u^2 - N(u)$, where $u>0$, and $N$
is as defined in \eqref{def-N}.  it follows from $(2.6)$ in
\cite{C83} and the inequality $ab \leq a^2 + b^2$ that
\begin{align*}
   \bigg| \int B(|X_\ve(s)|) d\xi  \bigg|
   \leq C |X_\ve|_2 |X_\ve|_{H^1}^{\frac{d}{d-2}}
   \leq C (|X_\ve|^2_2 + |X_\ve|_{H^1}^{\frac{2d}{d-2}}).
\end{align*}
Then
\begin{align*}
   \int N(|X_\ve|) d\xi
   \leq& \bigg| \int |X_\ve|^2 \log|X_\ve|^2 d\xi \bigg|
         + C (|X_\ve|^2_2 + |X_\ve|_{H^1}^{\frac{2d}{d-2}}).
\end{align*}
Hence by \eqref{esti-bdd-h1-gloabl} and Lemma \ref{Esti-Log},
\begin{align*}
    \bbe \sup\limits_{0\leq t\leq T}
    \bigg| \int N(|X_\ve|)d\xi \bigg|^p \leq C(T,p) <\9,
\end{align*}
which along with \eqref{V-N} implies \eqref{esti-v-p}, thereby
proving \eqref{esti-u-p}.

As regards \eqref{esti-u'-p}, we note that
\begin{align*}
   &\| e^W \calg_\ve(y_\ve)\|_{L^{p}(\Omega;C([0,T]; U'))} \nonumber \\
    \leq& 2|\lbb| \|X_\ve\|_{L^{p}(\Omega;C([0,T]; L^2))}
    + 2|\lbb| \| I_{\{|X_\ve| > e^{-3}\}} X_\ve L_\ve(X_\ve)\|_{L^{p}(\Omega;C([0,T]; L^2))}  \nonumber\\
    &+2|\lbb| \| I_{\{|X_\ve|\leq e^{-3}\}} X_\ve L_\ve(X_\ve) \|_{L^{p}(\Omega;C([0,T]; V'))}
\end{align*}
Since for $\xi\in\{|X_\ve| > e^{-3}\}$,
\begin{align} \label{xlve-leq1}
   |X_\ve(\xi) L_\ve(X_\ve (\xi)) |
   \leq |X_\ve(\xi)L(X_\ve(\xi)) |
   \leq  C_\delta (|X_\ve(\xi)| + |X_\ve(\xi)|^{1+\delta}),
\end{align}
by the Sobolev imbedding theorem and \eqref{esti-u-p}, it follows
that for $0\leq \delta\leq \frac{2}{d-2}$,
\begin{align*}
    &\|X_\ve\|_{L^{p}(\Omega;C([0,T]; L^2))}
    +  \| I_{\{|X_\ve|> e^{-3}\}} X_\ve L_\ve(X_\ve)\|_{L^{p}(\Omega;C([0,T];
    L^2))} \\
    \leq& C( \|X_\ve\|_{L^{p}(\Omega;C([0,T]; L^2))}
          + \|X_\ve\|^{1+\delta}_{L^{(1+\delta)p}(\Omega;C([0,T];
          L^{2(1+\delta)}))}) \\
    \leq& C( \|X_\ve\|_{L^{p}(\Omega;C([0,T]; L^2))}
          + \|X_\ve\|^{1+\delta}_{L^{(1+\delta)p}(\Omega;C([0,T];
          H^1))}) \leq C(T,p)<\9,
\end{align*}
where $C(T,p)$ is independent of $\ve$.

Moreover, since $\wt{N}$ is increasing, by Lemma \ref{Pro-Lve}
$(i)$, similarly to \eqref{wtN-V} we have
\begin{align*}
    \int \wt{N} (-2I_{\{|X_\ve|\leq e^{-3} \}} X_\ve L_\ve(X_\ve) ) d\xi
    \leq & \int I_{\{|X_\ve|\leq e^{-3} \}} \wt{N} (-X_\ve
    L(|X_\ve|^2)) d\xi \\
    \leq& 2 \max\{\| X_\ve\|_V, \| X_\ve\|^2_V \}.
\end{align*}
Then, as in \eqref{bdd-lp'-v'},
\begin{align*}
   &\| I_{\{|X_\ve|\leq e^{-3}\}} X_\ve L_\ve(X_\ve) \|_{L^{p}(\Omega;C([0,T];
   V'))} \nonumber \\
   \leq& C_T(\|X_\ve\|^2_{L^{2p}(\Omega;C([0,T]; V))} + 1)
   \leq C(T,p) < \9,
\end{align*}
where the last step is due to \eqref{esti-u-p}, thereby proving
\eqref{esti-u'-p}. The proof of Proposition \ref{Prop-yve-op} is now
complete. \hfill $\square$

\section{Proof of Theorem \ref{thm-y}} \label{Proof-Main}

Let us start with the lemma below.

\begin{lemma} \label{Esti-U'}
Let $L_{\ve}$ be defined as in \eqref{def-l} and $p \geq 3$. For any
$X\in L^p(\Omega \times(0,T); U)$,
\begin{align} \label{app-gev-g}
    \| XL_\ve(X) - X L(X)  \|_{L^{p'}(\Omega \times (0,T); U')} \to 0,\ \ as\ \ve \to
    0.
\end{align}
\end{lemma}

{\it \bf Proof.} First note that $X L_\ve(X) \to XL(X)$ a.e., as
$\ve \to 0$, and
\begin{align} \label{u'-split}
     &\| X L_\ve(X) - XL(X)
    \|_{L^{p'}(\Omega \times (0,T); U')} \nonumber \\
    \leq& \| I_{\{|X|\leq e^{-3}\}} (X L_\ve(X)  - XL(X))\|_{L^{p'}(\Omega \times (0,T); V')}  \nonumber \\
         &+ \| I_{\{|X| > e^{-3}\}} (X L_\ve(X)  - XL(X))
         \|_{L^{p'}(\Omega \times (0,T); H^{-1})}.
\end{align}

Since $|X L_\ve (X)| \leq |X L(X)|$, and as in the proof of
\eqref{bdd-lp'-h-1}, the Sobolev imbedding theorem and the
H\"{o}lder inequality imply that for $0<\delta\leq
\min\{\frac{2}{d-2} , p-2\}$,
\begin{align} \label{esti-w'-2}
   & \| I_{\{ |X|> e^{-3}\}} X L(X)
  \|_{L^{p'}(\Omega \times (0,T);
    H^{-1})} \nonumber  \\
   \leq&  C (\|X\|_{L^{p'}(\Omega \times (0,T);
    L^2)}
        + \|X\|^{1+\delta}_{L^{(1+\delta)p'}(\Omega \times(0,T);
        L^{2(1+\delta)})}) \nonumber \\
   \leq& C(T,p) (\|X\|_{L^p(\Omega\times (0,T);H^1)} + \|X\|^{1+\delta}_{L^p(\Omega\times (0,T);H^1)}
   )<\9,
\end{align}
the dominated convergence theorem implies that, as $\ve \to 0$,
\begin{align} \label{lve-l-h-1}
     &\| I_{\{|X| >e^{-3}\}} (X L_\ve(X) - X L(X))
         \|_{L^{p'}(\Omega\times(0,T); H^{-1})} \nonumber \\
     \leq& \| I_{\{|X| > e^{-3}\}} (X L_\ve(X) - X L(X))
         \|_{L^{p'}(\Omega\times(0,T); L^2)}   \to 0.
\end{align}

For the last term in the right hand side of \eqref{u'-split}, note
that since $\wt{N}$ is increasing, by Lemma \ref{Pro-Lve} $(i)$ and
\eqref{N'-N}
\begin{align} \label{nxlve}
   \wt{N}(-2|X(\xi)|L_\ve(X(\xi)) )
   \leq \wt{N} (-|X(\xi)|L(|X(\xi)|^2)) \leq 2 N(|X(\xi)|).
\end{align}
Moreover, by \eqref{V-N} and H\"{o}lder's inequality
\begin{align*}
    \| \int N(|X(\xi)|) d\xi \|_{L^{p'}(\Omega\times (0,T))}
    \leq& C_T (\|X\|^2_{L^p(\Omega\times(0,T);V)} +1)<\9,
\end{align*}
which implies that $N(|X|) \in L^1(\bbr^d)$, $\bbp \otimes dt$-a.e.
Then, it follows from the dominated convergence theorem that $\bbp
\otimes dt$-a.e.
\begin{align*}
     \int \wt{N}( -2X L_\ve(X) I_{\{|X|\leq e^{-3}\}})  d\xi
     \to \int \wt{N}(-2X L (X)I_{\{|X|\leq e^{-3}\}} )d\xi,
\end{align*}
which yields by \cite[(2.8)]{C83} that
\begin{align*}
     -X L_\ve(X) I_{\{|X|\leq e^{-3}\}}
     \to -X L (X)I_{\{|X|\leq e^{-3}\}},\ \ in\ V'.
\end{align*}
Since by Lemma \ref{Pro-Lve} $(i)$, $\|-XL_\ve(X)I_{\{|X|\leq
e^{-3}\}}\|_{V'} \leq  \|-XL(X)I_{\{|X|\leq e^{-3}\}}\|_{V'} $, and
as in \eqref{bdd-lp'-v'}, we have
$
    \|-X L (X)I_{\{|X|\leq
    e^{-3}\}}\|_{V'} \ \in L^{p'}(\Omega \times (0,T))
$. Again, we apply the dominated convergence theorem and get
\begin{align} \label{lve-l-v'}
    \| I_{\{|X|\leq e^{-3}\}} (X L_\ve(X)  - XL(X))\|_{L^{p'}(\Omega \times (0,T);
    V')}  \to 0,\ \ \ve\to 0.
\end{align}

Consequently, \eqref{app-gev-g} follows from \eqref{u'-split},
\eqref{lve-l-h-1} and \eqref{lve-l-v'}. The proof of Lemma
\ref{Esti-U'} is thus complete.
\hfill $\square$ \\

{\bf Proof of Theorem \ref{thm-y}.} For any $p\geq 2$, by the
uniform estimates \eqref{esti-u-p} and \eqref{esti-u'-p}, we have
along a subsequence $\{\ve_n\} \to 0$,
\begin{align}
    &e^Wy_{\ve_n}  \overset{\omega^*}{\rightharpoonup} e^W \wt{y},\ in\ L^p(\Omega; L^\9(0,T;U)),
    \label{app-weak.2-C} \\
    & e^W\calg_{\ve_n}(y_{\ve_n}) \overset{\omega^*}{\rightharpoonup} e^W\eta,\ in\  L^p(\Omega; L^\9(0,T;U')),
    \label{app-weak.1-C}
\end{align}
where  $\overset{\omega^*}{\rightharpoonup}$ stands for weak-star
convergence.

In particular, $e^W\wt{y} \in L^\9(0,T; U)$, and $e^W\eta \in
L^\9(0,T; U')$, $\bbp$-a.s. Since by Hypothesis $(H)$, for any $u\in
U$, we have $\|e^{-W}u \|_{U} \leq c(t) \|u\|_{U}$, where $c(t)=
\sqrt{2}(|e^{-W(t)}|_{L^\9} + |\na e^{-W(t)}|_{L^\9})$. It follows
that
\begin{align*}
\wt{y} \in L^\9(0,T; U),\ \ \  \eta \in L^\9(0,T; U'),\ \ \bbp-a.s.
\end{align*}

Moreover, for any $p\geq 3$, since $L^p(\Omega; L^\9(0,T;U)) \subset
L^p(\Omega \times (0,T);U))$ and $L^p(\Omega; L^\9(0,T;U')) \subset
L^{p'}(\Omega \times (0,T);U')$, we have (selecting a further
subsequence if necessary)
\begin{align}
    &y_{\ve_n}  \overset{\omega}{\rightharpoonup} \wt{y},\ in\ \mathcal{U},
    \label{app-weak.2} \\
    & \calg_{\ve_n}(y_{\ve_n}) \overset{\omega}{\rightharpoonup} \eta,\ in\  \mathcal{U}'.
    \label{app-weak.1}
\end{align}
where $\overset{\omega}{\rightharpoonup}$ means weak convergence.

We next take the limit in the approximating equation
\eqref{app-equa-y}. Set
$$F_n(y_{\ve_n})(t):= i e^{-W(t)} \Delta (e^{W(t)}y_{\ve_n}(t)) +
(4i\lbb |\lbb|t + \wh{\mu})y_{\ve_n}(t) +
\calg_{\ve_n}(y_{\ve_n}(t)),$$ and
$$F(\wt{y})(t):=i e^{-W(t)} \Delta (e^{W(t)}\wt{y}(t)) + (4i\lbb |\lbb|t +
\wh{\mu})\wt{y}(t) + \eta(t),$$ where $t\in [0,T]$. Then, from
\eqref{app-weak.2} and \eqref{app-weak.1} it follows that
$F_{n}(y_{\ve_n}) \overset{\omega}{\rightharpoonup} F(\wt{y})$, in
$\calu'$. Thus, for any $u\in U$, $\vf \in L^\9([0,T]\times
\Omega)$, by \eqref{app-equa-y}
\begin{align*}
   &\bbe \int_0^T\ {}_{U'}\<e^{W(t)}\wt{y}(t), \vf(t)u\>_{U} dt \\
   =& \lim\limits_{n \to \9} \bbe \int_0^T\ {}_{U'}\<e^{W(t)}\wt{y}_{\ve_n}(t), \vf(t)u\>_{U}
   dt\\
   =& \bbe \int_0^T\ \<e^{W(t)}x, \vf(t)u\>_H dt
      -  \lim\limits_{n\to \9} \bbe \int_0^T \int_0^t\ {}_{U'}\<F_n(y_{\ve_n})(s),
      \ol{e^{W(t)}}\vf(t)u\>_{U} ds dt.
\end{align*}
Note that, the second term in the right hand side above is equal to
\begin{align*}
&\lim\limits_{n\to \9} \bbe \int_0^T \
{}_{U'}\<e^{W(s)}F_n(y_{\ve_n})(s),
      \int_s^T\ol{e^{W(t)-W(s)}}\vf(t)u dt\>_{U} ds \\
=& \bbe \int_0^T \ {}_{U'}\<e^{W(s)}F(\wt{y})(s),
      \int_s^T\ol{e^{W(t)-W(s)}}\vf(t)u dt\>_{U} ds \\
=& \bbe \int_0^T {}_{U'}\<e^{W(t)}\int_0^t F(\wt{y})(s)ds,
\vf(t)u\>_U dt.
\end{align*}
It follows that for any $u\in U$, $\vf \in L^\9([0,T]\times
\Omega)$,
\begin{align*}
   \bbe \int_0^T\ {}_{U'}\<e^{W(t)}\wt{y}(t), \vf(t)u\>_{U} dt
   =\bbe \int_0^T\ {}_{U'}\<e^{W(t)}\(x-\int_0^t F(\wt{y})(s)ds\), \vf(t)u\>_{U}
   dt.
\end{align*}
Thus, $\wt{y} = x- \int_0^{\cdot} F(\wt{y})(s) ds,$ in $\calu'$.

Set $y(t):= x -\int_0^t F(\wt{y})(s)ds$, $t\in[0,T]$. Then, $y\in
AC([0,T]; U')$, $\bbp$-a.s., and $y=\wt{y}$ in $\calu'$, which
implies that $y = \wt{y}$, $\bbp \otimes dt$-a.e. Since for each
$t\in [0,T]$, $\int_0^t F(y)(s)ds = \int_0^t F(\wt{y})(s)ds$,
$\bbp$-a.s., by the continuity of $t \mapsto \int_0^t F(y)(s) -
F(\wt{y})(s) ds$, we thus have that $\bbp$-a.s. for all $t\in
[0,T]$, $\int_0^t F(y)(s)ds = \int_0^t F(\wt{y})(s)ds$, which yields
that $y (t) = \wt{y}(t)$, for all $t\in[0,T]$, $\bbp$-a.s.

Therefore,
\begin{align} \label{equa-y-eta}
    y(t)=x - \int_0^t (i e^{-W(s)}\Delta(e^{W(s)}y(s)) + &(4i \lbb|\lbb| s + \wh{\mu})y(s) +
    \eta(s))ds,\nonumber\\
    &\qquad for\ all\ t\in [0,T],\ \bbp-a.s.
\end{align}

Moreover, taking into account  \eqref{equa-y-eta},
\eqref{app-weak.2} and \eqref{app-weak.1}, we have $y\in
W^{1,p'}(0,T; U') \cap L^p([0,T]; U)$,  $\bbp$-a.s., which implies
that $y\in C([0,T]; H)$, $\bbp$-a.s.\\

In order to prove that $y$ is a solution to \eqref{equa-y*}, we need
to show that
\begin{align} \label{eta-gy}
     \eta = \calg(y).
\end{align}
For this purpose, it suffices to prove that
\begin{align} \label{mono-lim}
   \limsup\limits_{n\to \9} \int_0^T Re\  {}_{\mathcal{U}_t'} \<\calg_{\ve_n} (y_{\ve_n}),
    y_{\ve_n} \>_{\mathcal{U}_t} dt
    \leq  \int_0^T Re\ {}_{\mathcal{U}_t'} \<\eta,
y\>_{\mathcal{U}_t} dt,
\end{align}
where $\calu_t$ and $\calu'_t$ are defined as in \eqref{def-calU}
and \eqref{def-calU'} respectively, but with $T$ replaced by $t$.
Indeed, by the monotonicity of $\calg_{\ve_n}$, for any positive
function  $\varphi \in C([0,T])$,
\begin{align*}
   \int_0^T Re\ {}_{\calu'_t} \< \calg_{\ve_n}(y_{\ve_n})-\calg_{\ve_n}(u),
   y_{\ve_n}-u
\>_{\calu_t} \vf(t) dt \geq 0,\ \ u\in \calu.
\end{align*}
Then, it follows from \eqref{app-gev-g}, \eqref{app-weak.2} and
\eqref{app-weak.1} that
\begin{align*}
   \limsup\limits_{n\to\9}
     Re\ {}_{\calu'_t} \< \calg_{\ve_n}(y_{\ve_n})-\calg_{\ve_n}(u), y_{\ve_n}-u
\>_{\calu_t}
    \leq  Re\ {}_{\calu'_t} \< \eta-\calg(u), y-u
\>_{\calu_t}.
\end{align*}
Moreover,
\begin{align*}
   &| Re\ {}_{\calu'_t} \< \calg_{\ve_n}(y_{\ve_n})-\calg_{\ve_n}(u), y_{\ve_n}-u
\>_{\calu_t}| \nonumber \\
   \leq& \sup\limits_{n\geq1} (\|\calg_{\ve_n}(y_{\ve_n})\|_{\calu'} + \|\calg_{\ve_n}(y)\|_{\calu'})
    (\|y_{\ve_n}\|_{\calu'} + \|y\|_{\calu'}) < \9.
\end{align*}
Hence, by Fatou's lemma,
\begin{align*}
    0 \leq&  \limsup\limits_{n\to\9} \int_0^T Re\ {}_{\calu'_t} \< \calg_{\ve_n}(y_{\ve_n})-\calg_{\ve_n}(u),
    y_{\ve_n}-u
\>_{\calu_t} \vf(t) dt \\
    \leq&  \int_0^T \limsup\limits_{n\to\9} Re\ {}_{\calu'_t} \< \calg_{\ve_n}(y_{\ve_n})-\calg_{\ve_n}(u),
    y_{\ve_n}-u
\>_{\calu_t} \vf(t) dt \\
    =&  \int_0^T  Re\  {}_{\calu'_t} \< \eta-\calg(u), y-u
\>_{\calu_t} \vf(t) dt.
\end{align*}
As the integrand is continuous in $t$, and $\vf$ is an arbitrary
positive continuous function, we deduce that,
\begin{align*}
     Re\  {}_{\calu'} \< \eta-\calg(u), y-u
\>_{\calu} \geq 0,
\end{align*}
which implies \eqref{eta-gy} by the maximal monotonicity of $\calg$.
\\

For the proof of \eqref{mono-lim}, we note that by
\eqref{app-equa-y} we have, via It\^{o}'s formula,
\begin{align}  \label{mono-lim*}
   \int_0^T \bbe |e^{W(t)}y_{\ve_n}(t)|_{2}^2 dt
   =&|x|_{2}^2 T - 4|\lbb|  \int_0^T \bbe \int_0^t |y_{\ve_n}(s)|_2^2 ds dt
\nonumber \\
   =& |x|_{2}^2 T - 2 \int_0^T Re\ {}_{\calu'_t}\<
\calg_{\ve_n}(y_{\ve_n}),y_{\ve_n}\>_{\calu_t} dt.
\end{align}
Moreover, as in the proof of \cite[Lemma 8.1]{BR14}, applying
It\^o's formula to \eqref{equa-y-eta} we derive
\begin{align} \label{app.2}
    \int_0^T \bbe |e^{W(t)}y(t)|^2_{2} dt
   =&|x|^2_{2} T
      -2 \int_0^T Re\ {}_{\calu_t'}\<\eta,y\>_{\calu_t} dt.
\end{align}
Thus, by \eqref{app.2}, \eqref{app-weak.1} and \eqref{mono-lim*} we
derive that
\begin{align*}
   \int_0^T Re\ {}_{\calu_t'}\<\eta, y\>_{\calu_t} dt
   =&-\frac 12 \int_0^T \bbe |e^{W(t)} y(t)|_2^2 dt + \frac 12 |x|_2^2 T \\
   \geq& \limsup\limits_{n\to\9} \(-\frac 12 \int_0^T \bbe |e^{W(t)}y_{\ve_n}(t)|_2^2 dt + \frac 12
   |x|_2^2 T \) \\
   =&  \limsup\limits_{n\to\9} \int_0^T Re\ {}_{\calu_t'} \<\calg_{\ve_n}(y_{\ve_n}), y_{\ve_n}\>_{\calu_t} dt ,
\end{align*}
which yields \eqref{mono-lim}  as claimed, thereby proving
\eqref{eta-gy}.

Therefore, $y$ is a solution to \eqref{equa-y*} in the sense of
Definition \ref{def-y}. Moreover, the estimates
\eqref{Bdd-U-y}-\eqref{Bdd-U'-dy} follow immediately from
\eqref{app-weak.2-C}, \eqref{app-weak.1-C} and \eqref{equa-y-eta}.

It is left to prove the uniqueness, which follows from the
monotonicity.  In fact, given any two solutions $y_1,y_2$ to
\eqref{equa-y*}, setting $X_i = e^W y_i$, $i=1,2$, by the Ito
formula,  we obtain similar formula as in \eqref{difference} but
with $\ve=0$. Thus, it follows from \eqref{mono-lve-0} with $\ve=0$
and similar arguments as those below \eqref{difference} that
$X_1(t)= X_2(t)$, $\forall t \in[0,T]$, $\bbp$-a.s. The proof of
Theorem \ref{thm-y} is, therefore, complete. \hfill $\square$

\section{Appendix}

{\it \bf Proof of Lemma \ref{Pro-Lve}.} $(i)$. First note that, for
each $0<\ve <1$ fixed,
\begin{align*}
  \frac{d}{du} L_\ve(u) = \frac{1-\ve^2}{(\ve + u)(1+ \ve u)} \geq
  0,\ u>0,
\end{align*}
which implies that $L_\ve(u)$ is increasing with $u$, and so
$|L_\ve(u)| \leq |\log \ve|$.

Similarly, for each $u>0$ fixed,
\begin{align*}
  \frac{d}{d\ve} L_\ve(u) = \frac{1- u^2}{(\ve + u)(1+ \ve u)},\ u>0,
\end{align*}
which yields that $ \ve \mapsto L_\ve(u)$ is increasing with $\ve$
if $u\in[0,1]$, but decreasing if $u\in [1,\9)$. Hence, for $u\in
[0,1]$, we have $\log u \leq L_\ve(u) \leq 0$, and for $u\in
[1,\9)$, $0\leq L_\ve(u) \leq \log u$. Therefore, we obtain for all
$u>0$, $|uL_\ve(u)| \leq |u L(u)|$.

$(ii)$. We may assume $0<|u_2|\leq |u_1|$ without loss of
generality. Note that
\begin{align*}
     u_1 L_{\ve}(u_1) - u_2 L_{\ve}(u_2)
     = u_2
\(L_{\ve}(u_1) -L_{\ve}(u_2) \)
       +(u_1-u_2)L_{\ve}(u_1).
\end{align*}
Since $|L_{\ve}(u_1)|\leq |\log \ve|$, and
\begin{align} \label{Lve-diff-0}
   \left|L_{\ve}(u_1) - L_{\ve}(u_2) \right|
    \leq&  \frac{1+\ve|u_2|}{|u_2|+\ve}\ \left|\frac{|u_1|+\ve}{1+\ve|u_1|} -
    \frac{|u_2|+\ve}{1+\ve|u_2|} \right|  \nonumber \\
    =& \left|\frac{(1-\ve^2)(|u_1|-|u_2|)}{(|u_2|+\ve)(1+\ve |u_1|)} \right| \nonumber \\
    \leq& (1-\ve^2) |u_2|^{-1}|u_1 - u_2|,
\end{align}
we obtain immediately \eqref{Lve-diff}.

$(iii)$. We assume $0<|u_2|\leq |u_1|$ without loss of generality.
Note that
\begin{align*}
   Im (\ov{u_1} - \ov{u_2}) (u_1 L_\ve(u_1) - u_2 L_\ve(u_2))
   = (Im (\ov{u_1} u_2)) (L_\ve(u_1) - L_\ve(u_2)),
\end{align*}
and
\begin{align*}
   |Im(\ov{u_1} u_2)|
   = \big|\frac{u_2(\ov{u_1}-\ov{u_2}) + \ov{u_2}(u_2-u_1)}{2i} \big|
    \leq |u_2||u_1 - u_2|.
\end{align*}
Thus, taking into account \eqref{Lve-diff-0} we obtain
\eqref{mono-lve-0}. \hfill
$\square$ \\

{\it \bf Proof of Lemma \ref{Evolu-A}.} This lemma follows
essentially from \cite{D94, D96}. Using the notations in \cite{D96},
we reformulate \eqref{equa-U} in form
\begin{align*}
   (\partial_t + i\D + \sum\limits_{j=1}^d b^jD_j +
   c) y = f,
\end{align*}
meant in the weak sense, where $D_j = - i\partial_j$, $b^j =
-2\partial_j W$, and $c= i \sum\limits_{j=1}^d  (\partial_j W)^2 +
i\D W  + 2|\lbb| + 4i\lbb|\lbb|t + \wh{\mu} $.

Since for each  $1\leq m \leq n$, $e_m \in C_b^\9$ and
$\beta(\cdot)$ is continuous, $\bbp$-a.s., we have $b^j, c \in
C_{\omega}([0,T];\calb^\9)$, where $\calb^\9$ is as defined in
Hypothesis $(H)$, and $C_{\omega}([0,T];\calb^\9) =\{g\in C([0,T];
C^\9), \{g(t,\cdot)\}_{0\leq t\leq T}\ is\ uniformly\ bounded\ in\
\calb^\9. \}$

Moreover, under Hypothesis  $(H)$,
\begin{align*}
  |Re\ b(t,\xi)| \leq \(2\sum\limits_{m=1}^n |\mu_m|\sup\limits_{t\leq T}
  |\beta_m(t)|\) \lbb(|\xi|).
\end{align*}

Hence, the conditions in \cite[Theorem 1.1]{D96} are verified, and
we obtain the existence and uniqueness of the evolution operators
$U(t,s)$.

Furthermore, as remarked by the author in \cite{D96}, the results in
\cite{D94} holds also for the time-dependent coefficients. Thus,
similarly to \cite[(1.6)]{D94}, we have the estimates
\eqref{Esti-U-0} and \eqref{Esti-U}.

Finally, the measurabilities of the processes $U(\cdot, s)x$ and
$C_t$, $t\geq 0$, can be proved similarly as in the proof of Lemma
$3$ and Lemma $4$ in \cite{BRZ14} (see also \cite[Lemma 1.2.1, Lemma
1.2.3]{Z14}). The proof is now
complete. \hfill $\square$ \\

{\it \bf Proof of \eqref{Ito-Fn}.} Since the nonlinearity $X_\ve
L_\ve(X_\ve) \in L^2 \subset H^{-1}$, we can use similar arguments
as in the proof of \cite[Lemma 2.4, Proposition 6.1]{BRZ14.2} to
derive that $X_\ve := e^W y_\ve$ satisfies $\bbp$-a.s. for all
$t\in[0,T]$,
\begin{align}
    X_\ve(t) =& x -i \int_0^t \Delta X_\ve ds - 2\lbb i \int_0^tX_\ve L_\ve(X_\ve) ds \nonumber \\
   &- \int_0^t(4i \lbb |\lbb|s + 2|\lbb| +\mu) X_\ve ds + \sum\limits_{j=1}^n \int_0^t X_\ve \phi_j
   d\beta_j(s), \label{equa-x-le}
\end{align}
where the equation is taken in $H^{-1}$.

 Proceeding as in \cite{K10} and \cite{BRZ14.2}, we set
$h_\delta = h \ast \psi_\delta$ for any locally integrable function
$h$ mollified by $\psi_\delta$, where $\psi_\delta = \delta^{-d}
\psi(\frac{x}{\delta})$ and $\psi$ is a real-valued, nonnegative,
compactly supported smooth function with unit integral.

Taking convolution of both sides of \eqref{equa-x-le} with the
mollifiers $\psi_\delta$, we have for each $\xi \in \bbr^d$ that
\begin{align}
   (X_\ve(t))_\delta(\xi) =& -i \int_0^t \Delta X_{\ve,\delta}(\xi) ds - 2\lbb i \int_0^t (X_\ve L_\ve(X_\ve))_\delta(\xi) ds \nonumber   \\
   &-  \int_0^t (4i \lbb |\lbb|s + 2|\lbb|) X_{\ve,\delta}(\xi) ds
    - \int_0^t (\mu X_\ve)_\delta (\xi) ds \nonumber \\
   &+  \sum\limits_{j=1}^n\int_0^t (X_\ve \phi_j)_\delta(\xi)
   d\beta_j (s),\ \ t\in[0,T], \label{equa-x-le-delta}
\end{align}
where $X_{\ve, \delta}= (X_\ve)_\delta$, and \eqref{equa-x-le-delta}
holds on a set $\Omega_\xi\in \mathscr{F}$ with
$\bbp(\Omega_\xi)=1$.

Since for any locally integrable function $h$, $h_\delta(\xi)$ is
continuous in $\xi$, using the boundedness of the $H^1$-norm in
\eqref{esti-bdd-h1-gloabl} and similar arguments as in the proof of
\cite[Lemma 5.1]{BRZ14.2} and \cite[Lemma 2.3.11]{Z14}, we can prove
the continuity in $\xi$ of all terms in \eqref{equa-x-le-delta}.
Thus, \eqref{equa-x-le-delta} holds on a full probability set
$\wt{\Omega}\in \mathscr{F}$, which is independent of $\xi\in
\bbr^d$. For simplicity, below we omit the argument $\xi$ in
\eqref{equa-x-le-delta}.

Now, applying It\^{o}'s formula to the real valued function
$F_m(|X_{\ve,\delta}|^2)$, then integrating over $\bbr^d$,
interchanging the integrals and integrating by parts,  we obtain
\begin{align*}
  &\int F_m(X_{\ve, \delta}(t)) d\xi \\
  =& \int F_m(x) d\xi
     - 2\int_0^t  \int g_m(|X_{\ve, \delta}|^2)
        Re (\ov{X_{\ve,\delta}} \na X_{\ve,\delta})
        Im (\ov{X_{\ve,\delta}} \na X_{\ve,\delta}) d\xi ds\\
  &+ 4\lbb Im \int_0^t \int ( L_{1/m}(|X_{\ve,\delta}|^2) +1)
  \ov{X_{\ve,\delta}} (X_\ve L_\ve(X_\ve))_\delta d\xi ds\\
  & - 4|\lbb| \int_0^t \int
  (L_{1/m}(|X_{\ve,\delta}|^2)+1)|X_{\ve,\delta}|^2 d\xi ds\\
  & -2 \int_0^t \int (L_{1/m}(|X_{\ve,\delta}|^2)+1) Re (\ov{X_{\ve,\delta}} (\mu
  X_\ve)_\delta) d\xi ds\\
  & + \sum\limits_{j=1}^n \int_0^t \int  (L_{1/m}(|X_{\ve,\delta}|^2)+1) |(X_\ve
  \phi_j)_\delta|^2 d\xi ds\\
  & + \sum\limits_{j=1}^n  \int_0^t \int g_m(|X_{\ve, \delta}|^2) (Re (\ov{X_{\ve,
  \delta}}(X_{\ve}\phi_j)_\delta))^2 d\xi ds\\
  & + 2 \sum\limits_{j=1}^n  \int_0^t \int (L_{1/m}(|X_{\ve,\delta}|^2)+1) Re
  (\ov{X_{\ve,\delta}} (X_\ve \phi_j)_\delta) d\xi d\beta_j(s),
\end{align*}
where $ g_m(|X_{\ve, \delta}|^2) := 2(1-\frac{1}{m^2})(\frac{1}{m} +
|X_{\ve, \delta}|^2)^{-1} (1+ \frac{1}{m} |X_{\ve,
\delta}|^2)^{-1}$. (Note that, since $|L_{1/m}(|X_{\ve,\delta}|^2)|
\leq \log m$, we can use the (stochastic) Fubini theorem to
interchange the integrals.)

Therefore, since $|L_{1/m}(|X_{\ve,\delta}|^2)| \leq \log m$,
$|g_m(|X_{\ve,\delta}|^2)| \leq 2 |X_{\ve,\delta}|^{-2}$,  and
$h_\delta \to h$ in $L^q$, for any $h\in L^q$, $q>1$, using the
boundedness of the $H^1$-norm in \eqref{esti-bdd-h1-gloabl} and the
generalized dominated convergence theorem, we can take the limit
$\delta \to 0$ above and consequently obtain \eqref{Ito-Fn}. \hfill
$\square$

\end{document}